\documentclass[twoside,centertags                 
]{amsart}

\usepackage[cmtip,color,all,
]{xy}

\usepackage{amssymb}
\usepackage{graphicx}
\usepackage{mathrsfs}

\newtheorem{thm}[equation]{Theorem}
\newtheorem{cor}[equation]{Corollary}
\newtheorem{lem}[equation]{Lemma}
\newtheorem{prop}[equation]{Proposition}
\theoremstyle{definition}
\newtheorem{defn}[equation]{Definition}
\theoremstyle{remark}
\newtheorem{rem}[equation]{Remark}

\newcommand{\C}[1]{\mathscr{#1}}
\newcommand{\CC}[1]{\mathbb{#1}}
\def\op{\operatorname{op}}

\def\catata{\operatorname{cat}}

\def\env{\operatorname{e}}

\def\r{\rightarrow} 
\def\rr{\Rightarrow} 
\def\into{\rightarrowtail}
\def\onto{\twoheadrightarrow}
\def\To{\longrightarrow}

\newcommand{\id}[1]{\operatorname{id}_{#1}}

\def\hom{\operatorname{Hom}}
\def\ho{\operatorname{Ho}}

\def\st{\stackrel} 

\def\unit{\mathbf{1}} 

\def\colim{\mathop{\operatorname{colim}}}

\numberwithin{equation}{section}

\newdir{ >}{{}*!/-7pt/@{>}}
\newdir{> }{{}*!/10pt/@{>}}
\newdir{>> }{{}*!/10pt/@{>>}}

\def\dos{\mathbf{2}}

\newcommand{\inj}[1]{\operatorname{inj}(#1)}
\newcommand{\cof}[1]{\operatorname{cof}(#1)}
\newcommand{\cell}[1]{\operatorname{cell}(#1)}
\newcommand{\cellr}[1]{\operatorname{cell_{r}}(#1)}

\newcommand{\cat}[2]{\operatorname{Cat}_{#2}(#1)}

\renewcommand{\graph}[2]{\operatorname{Graph}_{#2}(#1)}

\newcommand{\alg}[1]{\operatorname{Mon}(#1)}
\newcommand{\mor}[1]{\operatorname{Mor}(#1)}
\newcommand{\lmod}[1]{\operatorname{Mod}(#1)}
\newcommand{\mon}{\alg}
\newcommand{\ob}{\operatorname{Ob}}

\begin{document}

\title[Dwyer--Kan homotopy theory of enriched categories]{Dwyer--Kan homotopy theory of enriched categories}%
\author{Fernando Muro}%
\address{Universidad de Sevilla,
Facultad de Matem\'aticas,
Departamento de \'Algebra,
Avda. Reina Mercedes s/n,
41012 Sevilla, Spain}
\email{fmuro@us.es}
\urladdr{http://personal.us.es/fmuro}


\begin{abstract}
We construct a model structure on the category of small categories enriched over a combinatorial closed symmetric monoidal model category satisfying the monoid axiom. Weak equivalences are Dwyer--Kan equivalences, i.e.~enriched functors which induce weak equivalences on morphism objects and equivalences of ordinary categories when we take sets of connected components on morphism objects.
\end{abstract}

\maketitle
\tableofcontents


\section{Introduction}

Several authors have studied the homotopy theory of small categories enriched in a closed symmetric monoidal model category $\C V$ in the sense of Hovey \cite[Definition 4.2.6]{hmc}, e.g.~Bergner \cite{mcssc} for $\C V$ the category of simplicial sets, Tabuada \cite{cmqdgcat,htsc} for $\C V$ the categories of chain complexes over a commutative ring and symmetric spectra with the stable model structure defined in \cite{se}, and Amrani \cite{msctc} for $\C V=\operatorname{Top}$ the category of weakly Hausdorff compactly generated topological
spaces. They construct model structures on the category $\cat{\C V}{}$ of small $\C V$-enriched categories and $\C V$-enriched functors. Lurie \cite[Proposition A.3.2.4]{htt} obtained the first result for a general $\C V$ satisfying a set of conditions. His hypotheses, however, are somewhat restrictive and do not cover all previously mentioned examples, just the categories of simplicial sets and chain complexes over a field.

In all the previous approaches, the authors concentrate in a fixed class of weak equivalences in $\cat{\C V}{}$, commonly known as \emph{Dwyer--Kan equivalences} after \cite{slc}, which generalizes the notion of equivalence of categories, undestood as a fully-faithful and essentially surjective functor, in a straightforward way. 

A $\C V$-enriched functor $\varphi\colon \CC H\r \CC K$ between small $\C V$-enriched categories is \emph{homotopically fully faithful} if $\varphi(x,y)\colon \CC H(x,y)\r \CC K(\varphi(x),\varphi(y))$ is a weak equivalence in $\C V$ for any two objects $x,y\in \ob\CC H$. 

The homotopical notion of essentially surjective is slightly more complicated. Let $\unit$ be the tensor unit in $\C V$. The \emph{`connected components'} functor
$$\pi_{0}=\ho(\C V) (\unit,-)\colon \C V\To \operatorname{Set}$$
is lax symmetric monoidal, hence any $\C V$-category $\CC K$ has a \emph{`category of connected components'} $\pi_{0}\CC K$, obtained by applying $\pi_{0}$ to morphism objects. A $\C V$-enriched functor $\varphi\colon \CC H\r \CC K$ is \emph{homotopically essentially surjective} if $\pi_{0}\varphi\colon \pi_{0}\CC H\r \pi_{0}\CC K$ is essentially surjective in the usual sense.

A $\C V$-enriched functor $\varphi\colon \CC H\r \CC K$ is a Dwyer--Kan equivalence, or simply a \emph{DK-equivalence}, if it is homotopically fully faithful and homotopically essentially surjective. In particular, $\pi_{0}\varphi\colon \pi_{0}\CC H\r \pi_{0}\CC K$ is an equivalence of categories.

Weak equivalences are not enough to \emph{define} a model structure. It is also necessary to choose the fibrations or cofibrations. It is even enough to choose the trivial fibrations or trivial cofibrations, as long as they are included in the chosen class of weak equivalences. If one of these four classes of maps is chosen together with the weak equivalences, then there is a \emph{canonical} choice of fibrations and cofibrations, imposed by the lifting axiom, which may or may not satisfy the axioms of a model category. If they satisfy the axioms, we say that the choice \emph{gives rise to} a model structure.

In all previously mentioned cases, \emph{trivial fibrations} in  $\cat{\C V}{}$ are the   $\C V$-enriched functors $\varphi\colon \CC H\r \CC K$ which are surjective on objects and such that the morphism $\varphi(x,y)\colon \CC H(x,y)\r \CC K(\varphi(x),\varphi(y))$ is a trivial fibration in $\C V$ for all $x,y\in \ob\CC H$. We say that  $\cat{\C V}{}$ admits \emph{the Dwyer--Kan model structure} if these trivial fibrations and the DK-equivalences give rise to a model structure on  $\cat{\C V}{}$.

The following theorem is the main result of this paper.

\begin{thm}\label{main}
Let $\C V$ be a combinatorial closed symmetric monoidal model category satisfying  Schwede--Shipley's monoid axiom \cite[Definitions  3.3]{ammmc}. Then $\cat{\C V}{}$ admits the Dwyer--Kan model structure. Moreover, this model structure is combinatorial. 
\end{thm}

Berger and Moerdijk \cite{htec} approach the problem of endowing $\cat{\C V}{}$ with a model structure in a different way. 
Trivial fibrations and fibrant objects do not \emph{define} a model structure. They yield a canonical choice for cofibrations by the lifting axiom, but not for fibrations or weak equivalences. This is a priori an inconvenience since, after all, weak equivalences determine the whole homotopy theory. Nevertheless, a result of Joyal ensures that if two model structures have the same  trivial fibrations and  fibrant objects then they must coincide \cite[Proposition E.1.10]{tqca}.

Berger and Moerdijk say that $\cat{\C V}{}$ admits  \emph{the canonical} model structure if there exists a model structure on it with the same trivial fibrations as above and such that $\CC H$ is fibrant if and only if $\CC H(x,y)$ is fibrant in $\C V$ for all $x,y\in \ob \CC H$. These are precisely the fibrant objects in all previously known cases. They find hypotheses on $\C V$ under which $\cat{\C V}{}$ admits the canonical model structure and, under the same hypotheses, show through a series of non-trivial results that the weak equivalences of this model structure are the DK-equivalences. In particular, the canonical model structure coincides with the Dwyer--Kan model structure in these cases.

Their hypotheses are more general than Lurie's,  in particular they cover Tabua\-da's and Amrani's results. However, they are still somewhat restrictive, they do not include other examples of interest such as diagram spectra with the positive stable model structure \cite{mcds} or, more generally, the underlying model category $\C V$ of a homotopical algebra context in the sense of To\"en and Vezzosi \cite{hagII}. Therefore Berger--Moerdijk's theorem cannot be used to define moduli spaces of enriched categories in general homotopical algebraic geometry contexts. The precise conditions, of technical nature,  include in an essential way right properness, cofibrancy of the tensor unit $\unit$, and a kind of compact generation which is not preserved under localizations. The rest of conditions are satisfied under the assumptions of Theorem \ref{main}.

The hypotheses of Theorem \ref{main} are satisfied by diagram spectra with the positive stable model structure, by the underlying category $\C V$ of any homotopical algebra context satisfying the monoid axiom, by any $\C V$ satisfying the assumptions of Lurie's proposition, and by all previously known specific examples, except for $\C V=\operatorname{Top}$, because it is not combinatorial. 

We do not know if, under the hypotheses of Theorem \ref{main}, the Dwyer--Kan model structure on $\cat{\C V}{}$ is always canonical. Fibrant $\C V$-categories $\CC H$ have fibrant morphism objects $\CC H(x,y)$ in $\C V$, however they must also satisfy a kind of isomorphism lifting condition which may not be redundant in some cases.  Berger--Moerdijk's conditions are sufficient, but a recent result of Stanculescu \cite{cmcpfo} shows that they are not necessary. He shows the existence and the coincidence of the Dwyer--Kan and the canonical model structures under certain hypotheses on $\C V$, which are in general stronger than  Berger--Moerdijk's (for instance, he asks $\C V$ to be simplicial), except for the fact that right properness is not required. Stanculescu's hypotheses are strictly stronger than ours.

We would like to mention a recent result of Haugseng \cite{reic} showing that, under some extra hypotheses, the $\infty$-category obtained by localizing $\cat{\C V}{}$ at DK-equivalences (in the simplicial way) is equivalent to the $\infty$-category of $\infty$-categories enriched in the localization of $\C V$ at weak equivalences. Theorem \ref{main} shows that this localization of $\cat{\C V}{}$ is always the localization of a combinatorial model category at weak equivalences.

In addition to our main theorem, we prove two other results showing that the Dwyer--Kan model structure is rather well behaved. The first one is a \emph{left properness} flavoured theorem.

\begin{thm}\label{lp}
Assume the  hypotheses of Theorem \ref{main} and that $\C V$ satisfies the strong unit axiom introduced in \cite[Definition A.9]{htnso2}. If
$$\xymatrix{
\CC H\ar@{>->}[r]^\varphi\ar[d]^\sim_\phi\ar@{}[rd]|{\text{push}}&\CC K\ar[d]^{\phi'}\\
\CC H'\ar@{>->}[r]_{\varphi'}&\CC K'
}$$
is a push-out in  $\cat{\C V}{}$
where $\CC H$ and $\CC H'$ have cofibrant morphism objects, $\varphi$ is a cofibration, and $\phi$ is a DK-equivalence, then $\phi'$ is also a DK-equivalence.
\end{thm}

The strong unit axiom is satisfied  in all examples we know. It holds if the tensor unit is cofibrant, but also if the tensor product of a cofibrant object and a weak equivalence is always a weak equivalence. Moreover, the hypotheses of Theorem \ref{lp} can be relaxed, see Theorem \ref{lp2}. It is enough to assume that $\CC H$ and $\CC H'$ have \emph{pseudo-cofibrant} morphism objects in the sense of \cite[Definition A.9]{htnso2}. An object $X$ in $\C V$ is pseudo-cofibrant if the functor $X\otimes -$ preserves cofibrations. Cofibrant objects are pseudo-cofibrant, but there are important examples of pseudo-cofibrant objects which need not be cofibrant, e.g.~the tensor unit. The morphism objects of any cofibrant $\C V$-category $\CC H$ are pseudo-cofibrant, but in general they are not cofibrant if the tensor unit is not cofibrant.

\begin{cor}\label{lp3}
Under the assumptions of Theorem \ref{main}, if all objects in $\C V$ are cofibrant then $\cat{\C V}{}$ is left proper.
\end{cor}

The second one concerns the change of base category. Recall that any lax symmetric monoidal functor $G\colon \C W\r\C V$ induces a functor $G\colon \cat{\C W}{}\r\cat{\C V}{}$ defined as $G$ on morphism objects.

\begin{thm}\label{transfer}
If $\C V$ and $\C W$ satisfy the assumptions of Theorem \ref{lp}  and $F\colon \C V\rightleftarrows\C W\colon G$ is a weak symmetric monoidal Quillen pair in the sense of \cite[Definition 3.6]{emmc} satisfying the pseudo-cofibrant axiom and the $\unit$-cofibrant axiom in \cite[Definition B.6]{htnso2} , then there is an induced Quillen pair $F^{\catata}\colon \cat{\C V}{}\rightleftarrows\cat{\C W}{}\colon G$ which is a Quillen equivalence if $F\colon\C V\rightleftarrows\C W\colon G$ is.
\end{thm}

The left adjoint $F^{\catata}$ is not  given by $F$ on morphism objects, unless $F$ is \emph{strong} symmetric monoidal. In general $F$ is just \emph{colax} symmetric monoidal. Nevertheless we have the following result. 

\begin{prop}\label{putifar}
Under the assumptions of Theorem \ref{transfer}, given a $\C V$-category $\CC H$
there is a natural isomorphism in $\ho\C W$
$$\mathcal LF^{\catata}(\CC H)(x,y)\cong \mathcal LF(\CC H(x,y))$$
for any $x,y\in\ob\CC H$.
\end{prop}

Here $\mathcal L$ denotes the total left derived functor of a left Quillen functor.

Weak monoidal Quillen adjunctions were introduced by Schwede and Shipley in \cite{emmc} in oder to obtain a Quillen equivalence between simplicial algebras and connective differential graded algebras out of the Dold--Kan equivalence $\operatorname{Mod}(\Bbbk)^{\Delta^{\op}}\rightleftarrows\operatorname{Ch}(\Bbbk)_{\geq 0}$, where $\Bbbk$ is a commutative ring. They were later used by Shipley \cite{hzas} to establish a zig-zag of Quillen equivalences between differential graded $\Bbbk$-algebras and algebras over the Eilenberg--MacLane ring spectrum $H\Bbbk$. Hence, the following results of Tabuada \cite{dgvsc,gscthhtm} can be recovered from Theorem \ref{transfer} in a unified way.

\begin{cor}
There is a Quillen equivalence $\cat{\operatorname{Mod}(\Bbbk)^{\Delta^{\op}}}{}\rightleftarrows\cat{\operatorname{Ch}(\Bbbk)_{\geq 0}}{}$ between the category of simplicial $\Bbbk$-linear categories and the category of connective $\Bbbk$-linear DG-categories.
\end{cor}

\begin{cor}
There is a zig-zag of Quillen equivalences between the category $\cat{\operatorname{Ch}(\Bbbk)}{}$ of $\Bbbk$-linear DG-categories and the category $\cat{\operatorname{Mod}(H\Bbbk)}{}$ of $H\Bbbk$-spectral categories 
\end{cor}



Since our principal results have been clearly stated in this introduction, the paper is written in a linear order, i.e.~we never invoke later results. For this reason, the three theorems above are proven in the three last sections of the paper. The proof of Theorem \ref{main} is very technical since we have to deal with push-outs in $\cat{\C V}{}$. It is precisely this what advises against cross-references. The other two theorems above do not make sense before Theorem \ref{main} is established. This is why they appear so late. Nevertheless, the Dwyer--Kan model structure of $\cat{\C V}{}$ is completely described in Section \ref{gtc}, which is the first one after two introductory sections. These introductory sections show how $\cat{\C V}{}$  is a bifibered category over sets and recall the known homotopy theory of the fibers $\cat{\C V}{S}$.

Throughout this paper, $\C V$ will denote a closed symmetric monoidal model category in the sense of Hovey \cite[Definition 4.2.6]{hmc}. In particular, it satisfies the push-out product axiom in  \cite[Definitions 3.1]{ammmc}, see also Definition \ref{ppax}  below. We will also assume that $\C V$  satisfies the monoid axiom, see \cite[Definition 3.3]{ammmc} or Definition \ref{monax}. 
In addition, we suppose that $\C V$ is combinatorial, i.e.~locally presentable. Most of our results also hold under weaker set-theoretical hypotheses, but we do not have relevant applications to non-combinatorial examples. Moreover, assuming once and for all that $\C V$ is combinatorial simplifies several arguments. We denote by $I$ and $J$ two fixed sets of generating cofibrations and  generating trivial cofibrations  of $\C V$, respectively. The symbols $\otimes$, $\unit$, and $\varnothing$ stand for the tensor product, the tensor unit, and the initial object of $\C V$, respectively.

In the purely categorical Sections \ref{fibered} and \ref{pushpush}, it would be enough that $\C V$ were a plain (co)complete closed symmetric monoidal category. 

\subsection*{Acknowledgements}

The author is grateful to Clemens Berger, Giovanni Caviglia, Ieke Moerdijk, Ji\v{r}\'i Rosick\'y and 
Daniel Sch\"appi
 for valuable discussions and comments conerning the contents of this paper. He 
was partially supported by the Andalusian Ministry of Economy, Innovation and Science under the grant FQM-5713,
by the Spanish Ministry of Education and
Science under the MEC-FEDER grant  MTM2010-15831, and by the Government of Catalonia under the grant SGR-119-2009.

\section{The bifibered categories of enriched categories and graphs}\label{fibered}

Small $\C V$-enriched categories and $\C V$-enriched functors between them will simply be called \emph{$\C V$-categories} and \emph{$\C V$-functors}, respectively. We refer the reader to \cite{borceux2,bcect} for the basics on (enriched) category theory.

The \emph{`set of objects' functor}
$$\ob\colon\cat{\C V}{}\To\operatorname{Set}$$
is a \emph{(Grothendieck) fibration} in the sense of \cite[Definition 8.1.3]{borceux2}. Some nice references for this topic in category theory are the corresponding chapters in \cite{borceux2,cna,elephant1}. The fiber at a given set $S$ is the subcategory $\cat{\C V}{S}\subset \cat{\C V}{}$ formed by the $\C V$-categories with object set $S$ and the $\C V$-functors which are the identity on objects. A map of sets $f\colon S\r S'$ induces a functor
\begin{equation}\label{*arriba}
f^{*}\colon\cat{\C V}{S'}\To \cat{\C V}{S} 
\end{equation}
characterized by the existence of a natural $\C V$-functor 
\begin{equation}\label{cartesian}
f^{*}\CC K\To \CC K, 
\end{equation}
actually a \emph{cartesian morphism} in $\cat{\C V}{}$, see  \cite[Definition 8.1.2]{borceux2}, given by $f$ on objects, which is the identity on morphisms
\begin{align*}
(f^{*}\CC K)(x,y)=\CC K(f(x),f(y)),\qquad x,y\in S.
\end{align*}
The functors $f^{*}$ define a \emph{splitting} of $\ob$, see  \cite[Definition 8.3.3]{borceux2}. 

We now want to show that this Grothendieck fibration is actually a \emph{bifibration}. This means that $\ob$, regarded as a functor between opposite categories $\ob\colon\cat{\C V}{}^{\op}\r\operatorname{Set}^{\op}$, is also a fibration. It is enough to prove that the functor $f^*$ in \eqref{*arriba} has a left adjoint, see  \cite[Lemma B1.4.5 (a)]{elephant1}. With this purpose, and for later use, we recall how $\cat{\C V}{S}$ arises from graphs.

\begin{defn}
A \emph{$\C V$-graph} $M$ is a set $\ob M$ together with a collection of objects $M=\{M(x,y)\}_{x,y\in \ob M}$ in $\C V$. A \emph{morphism of $\C V$-graphs} $\varphi\colon M\r N$ is a map $\varphi\colon\ob M\r\ob N$ together with a collection of morphisms  $$\{\varphi(x,y)\colon M(x,y)\To N(\varphi(x),\varphi(y))\}_{x,y\in\ob M}$$ in $\C V$. The category of $\C V$-graphs will be denoted by $\graph{\C V}{}$. 
\end{defn}

Notice that $\C V$-graphs are like $\C V$-categories but without composition or identities, and their morphisms are like $\C V$-functors. We will still refer to $M(x,y)$ as the \emph{morphis objects} of $M$, $x,y\in\ob M$, and say that $\varphi\colon M\r N$  is \emph{fully faithful} if $\varphi(x,y)\colon M(x,y)\r N(\varphi(x),\varphi(y))$ is an isomorphism for all $x,y\in\ob M$.

The `set of objects' functor for $\C V$-graphs
$$\ob\colon\graph{\C V}{}\To\operatorname{Set}$$
is clearly a Grothendieck bifibration. The fiber at a given set $S$ is the subcategory $\graph{\C V}{S}\subset \graph{\C V}{}$ of $\C V$-graphs with object set $S$ and morphisms which are the identity on objects. Notice that
$\graph{\C V}{S}$ 
is simply a product of $S\times S$ copies of $\C V$. A map of sets $f\colon S\r S'$ induces a functor
$$f^{*}\colon\graph{\C V}{S'}\To \graph{\C V}{S}$$
and cartesian morphisms $f^{*}M\r M$ defined as \eqref{*arriba} and \eqref{cartesian} above, and $f^*$ has a left adjoint
$$f_{*}\colon\graph{\C V}{S}\To \graph{\C V}{S'}$$
given by
$$f_*(M)(x,y)=\coprod_{
\begin{array}{c}
\\[-15pt]
\scriptstyle f(x')=x\\[-5pt]
\scriptstyle f(y')=y
\end{array}
} M(x',y').$$

\begin{defn}\label{otimess}
The  \emph{tensor product} of two $\C V$-graphs $M$ and $N$ with object set $S$ is defined by 
$$(M\otimes_S N)(x,y)= \coprod_{z\in S}M(z,y)\otimes N(x,z).$$
This tensor product endows $\graph{\C V}{S}$ with a (non-symmetric, in general) biclosed monoidal structure with \emph{tensor unit} 
$$\unit_{S}(x,y)=\left\{
\begin{array}{ll}
\unit,&x=y;\\
\varnothing,&x\neq y.
\end{array}
\right.$$
Associativity and unit constraints in $\graph{\C V}{S}$ are defined by the corresponding constraints in $\C V$.
\end{defn}

Clearly, $\cat{\C V}{S}$ is the category of monoids in $\graph{\C V}{S}$. Alternatively, it is the category of algebras over the monad associated to the adjoint pair
\begin{equation}\label{frees}
\xymatrix{\graph{\C V}{S}\ar@<.5ex>[r]^-{T_S}&\cat{\C V}{S},\ar@<.5ex>[l]^-{\text{forget}}} 
\end{equation}
where the right adjoint sends a $\C V$-category to its underlying $\C V$-graph, and the left adjoint $T_{S}$ is the \emph{free $\C V$-category} functor, recalled below in Section \ref{pushpush} in a more general context. 
Hence, the following result is a consequence of \cite[Theorem 4.5.6]{borceux2}.

\begin{prop}\label{*abajo}
The functor $f^{*}\colon \cat{\C V}{S'}\r \cat{\C V}{S}$ in \eqref{*arriba} has a left adjoint 
$$f_{*}\colon\cat{\C V}{S}\To \cat{\C V}{S'}.$$  
\end{prop}

The categories $\graph{\C V}{}$ and $\cat{\C V}{}$ are locally presentable, see \cite{vcatlplb}. The functors $\ob\colon\cat{\C V}{}\r\operatorname{Set}$ and $\ob\colon\graph{\C V}{}\r\operatorname{Set}$ preserve limits and colimits since they are both left and right adjoints, compare \cite[\S5]{vcatlplb}.  Moreover, the (co)limit of a diagram in $\graph{\C V}{}$  can be computed by moving the diagram to an appropriate fiber $\graph{\C V}{S}$, where $S$ is the (co)limit of the underlying diagram of object sets, and then computing the (co)limit in $\graph{\C V}{S}$. Similarly for $\cat{\C V}{}$. All (co)limits in $\graph{\C V}{S}$ are computed as in $\C V$ on morphisms. Limits and filtered colimits in $\cat{\C V}{S}$ are computed as in $\graph{\C V}{S}$, but other colimits, such as push-outs, are much harder to calculate, compare Section \ref{pushpush}.

In general, $f_{*}$ is also complicated to compute. However, if $f\colon S\hookrightarrow S'$ is injective, $f_{*}$ is surprisingly easy to describe: there is a natural  $\C V$-functor 
\begin{equation}\label{cocartesian}
\CC K\To f_{*}\CC K, 
\end{equation}
actually a \emph{cocartesian morphism} in $\cat{\C V}{}$,  given by $f$ on objects, which is the identity on morphisms, 
\begin{align*}
\CC K(x,y)=(f_{*}\CC K)(f(x),f(y)),\qquad x,y\in S,
\end{align*} 
and moreover, for $(x,y)\in (S'\times S')\setminus(f(S)\times f(S))$,
$$(f_{*}\CC K)(x,y)=\left\{
\begin{array}{ll}
\varnothing,&x\neq y,\\
\unit,&x=y.
\end{array}
\right.$$
The unit $\CC K\r f^{*}f_{*}\CC K$ is the identity and the counit $f_{*}f^{*}\CC K\r\CC K$ in the only $\C V$-functor which is the identy on morphism sets with source and target in $f(S)$.

%

Suppose $\C W$ is another category with the same structure and satisfying the same properties as $\C V$.
Any functor $G\colon\C W\r\C V$ induces a cartesian functor in the sense of \cite[Definition 8.2.1]{borceux2}
$$G\colon\graph{\C W}{}\To\graph{\C V}{}$$
which is defined as $G$ on morphism objects. This functor restricts on fibers
$$G_S\colon\graph{\C W}{S}\To\graph{\C V}{S}$$
to a product of $S\times S$ copies of $G$. Moreover, if $F\colon\C V\rightleftarrows\C W\colon G$  is an adjoint pair then the pairs
\begin{equation}\label{ap3}
\xymatrix{\graph{\C V}{}\ar@<.5ex>[r]^-{F}&\graph{\C W}{},\ar@<.5ex>[l]^-{G}}\qquad
\xymatrix{\graph{\C V}{S}\ar@<.5ex>[r]^-{F_S}&\graph{\C W}{S},\ar@<.5ex>[l]^-{G_S}} 
\end{equation}
are also adjoint.

If $G$ is lax symmetric monoidal then $G_S$ is lax monoidal for all $S$, 
so it lifts to monoids in a canonical way 
$$G_S\colon\cat{\C W}{S}\To\cat{\C V}{S}.$$
These functors assemble to a cartesian functor
\begin{equation}\label{kaker}
G\colon\cat{\C W}{}\To\cat{\C V}{}.
\end{equation}
However, if $F\colon\C V\rightleftarrows\C W\colon G$ is an adjunction, the lax symmetric monoidal structure of $G_S$ induces a \emph{colax} symmetric monoidal structure on $F_S$, which cannot be used to directly define a functor between categories of enriched categories in the opposite direction, unless this colax structure on $F_S$ happens to be strong. Nevertheless, by \cite[Theorem 4.5.6]{borceux2} there exists a Quillen pair
\begin{equation}\label{ap1}
\xymatrix{\cat{\C V}{S}\ar@<.5ex>[r]^-{F_S^{\catata}}&\cat{\C W}{S}.\ar@<.5ex>[l]^-{G_S}} 
\end{equation}
The left adjoint $F^{\catata}_S$ is misterious a priori, unless $F$ is strong, in which case $F^{\catata}_S=F_S$. By \cite[Corollaire 1.11.3]{cna}, the functors $F^{\catata}_S$ assemble to a left adjoint of $G$,
\begin{equation}\label{ap11}
\xymatrix{\cat{\C V}{}\ar@<.5ex>[r]^-{F^{\catata}}&\cat{\C W}{}.\ar@<.5ex>[l]^-{G}} 
\end{equation}

\section{Homotopy theory of enriched categories with fixed set of objects}\label{fixed}

The homotopy theory of the fibers $\cat{\C V}{S}$ has been considered in several sources, but it is still a very complicated problem to assemble these homotopy theories to a homotopy theory of the whole $\cat{\C V}{}$, due to the flabby nature of homotopical essential surjectivity. Nonetheless, in doing so, we will need some facts about $\cat{\C V}{S}$ that we recall in this section.

The category $\graph{\C V}{S}$ is a product of $S\times S$ copies of $\C V$, so it inherits a cofibrantly generated model structure with weak equivalences and (co)fibrations defined componentwise as in $\C V$. In order to exhibit sets of generating (trivial) cofibrations, given $x,y\in S$ we consider the adjoint pair
\begin{equation}\label{dosob}
\xymatrix{\C V\ar@<.5ex>[r]^-{(-)_{xy}}&\graph{\C V}{S}\ar@<.5ex>[l]^-{\text{forget}}} 
\end{equation}
where the forgetful functor is defined by $M\mapsto M(x,y)$. Its left adjoint sends an object $V$ in $\C V$ to the $\C V$-graph $V_{xy}$ with
$$V_{xy}(x',y')=\left\{
\begin{array}{ll}
V,&(x,y)=(x',y');\\
\varnothing,&(x,y)\neq(x',y').
\end{array}
\right.$$
The sets $$I_S=\bigcup_{x,y\in S}I_{xy}\qquad\text{and}\qquad  J_S=\bigcup_{x,y\in S}J_{xy}$$ are sets of generating cofibrations and generating trivial cofibrations in $\graph{\C V}{S}$, respectively.

The category $\cat{\C V}{S}$, being the category of monoids in the monoidal model category $\graph{\C V}{S}$, also inherits a cofibrantly generated model structure with weak equivalences and fibrations defined as in $\graph{\C V}{S}$. Such a model structure appears in \cite[Proposition 6.3]{emmc} under slightly stronger assumptions, and in \cite[Corollary 10.4]{htnso} under our assumptions. It is transferred along \eqref{frees}, which becomes a Quillen adjunction. The sets $$T_S(I_S)\qquad\text{and}\qquad T_S(J_S)$$  are sets of generating cofibrations and generating trivial cofibrations in $\cat{\C V}{S}$, respectively. The category $T_S(V_{xy})$ looks like the graph $V_{xy}$ if $x\neq y$, except that $T_S(V_{xy})(z,z)=\unit$ for all $z\in S$. If $x=y$, $T_S(V_{xx})(x,x)$ is the free monoid in $\C V$ on $V_{xx}$, $T_S(V_{xx})(z,z)=\unit$ if $z\neq x$, and $T_S(V_{xx})(z,z')=\varnothing$ for arbitrary $z\neq z'$. 

Clearly, the functor $f^*$ in \eqref{*arriba} induced by a map of sets preserves weak equivalences and fibrations, hence the adjunction $f_*\dashv f^*$ from Proposition \ref{*abajo} is a Quillen pair, and the same for $\C V$-graphs.

Let $\C W$ be, as in the previous section, a category with the same structure and satisfying the same properties as $\C V$. 
A Quillen pair $F\colon\C V\rightleftarrows\C W\colon G$ 
where the right adjoint $G$ is a lax symmetric monoidal functor 
induces Quillen pairs \eqref{ap3} and \eqref{ap1}. The misterious $F^{\catata}_S$ is related to $F_S$, and hence to $F$, via a natural transformation
\begin{equation}\label{chi}
\xymatrix{\cat{\C V}{S}\ar[r]^-{F_S^{\catata}}\ar[d]_{\text{forget}}&\cat{\C W}{S}\ar[d]^{\text{forget}}\\
\graph{\C V}{S}\ar[r]_-{F_S}&\graph{\C W}{S}\ar@{<=}(15,-4);(15,-10)_-{\chi}}
\end{equation}
which is the mate of the identity natural transformation in the following commutative square
$$\xymatrix{\cat{\C V}{S}\ar@{<-}[r]^-{G_S}\ar[d]_{\text{forget}}&\cat{\C W}{S}\ar[d]^{\text{forget}}\\
\graph{\C V}{S}\ar@{<-}[r]_-{G_S}&\graph{\C W}{S}}$$
The following result asserts that $\chi$ is almost a natural weak equivalence.

\begin{prop}\label{SS}
Suppose $F\colon\C V\rightleftarrows\C W\colon G$  is a weak symmetric monoidal Quillen pair in the sense of \cite[Definition 3.6]{emmc} satisfying the pseudo-cofibrant axiom and the $\unit$-cofibrant axiom in \cite[Definition B.6]{htnso2}. Assume also that  $\C V$ and $\C W$ satisfy the strong unit axiom in \cite[Definition A.9]{htnso2}. Then, if $\CC H$ is cofibrant in $\cat{\C V}{S}$, $\chi_{\CC H}\colon F_S(\CC H)\r F_S^{\catata}(\CC H)$ is a weak equivalence in $\graph{\C W}{S}$. 
\end{prop}

This proposition appears as \cite[Proposition 6.4 (1)]{emmc} under the extra assumption that tensor units are cofibrant in $\C V$ and $\C W$. Our hypotheses replace this cofibrancy assumption, see the appendices of \cite{htnso2}.

The Quillen pairs in \eqref{ap3} are Quillen equivalences as long as $F\colon\C V\rightleftarrows\C W\colon G$ is. 

\begin{prop}\label{qeq}
Under the hypotheses of Proposition \ref{SS}, if $F\colon\C V\rightleftarrows\C W\colon G$ is a Quillen equivalence then so is $F_S^{\catata}\dashv G_{S}$ in \eqref{ap1}.
\end{prop}

This follows easily from the previous proposition, see also \cite[Proposition 1.10 and Corollary E.5]{htnso2} for a more general approach.

\section{Generating (trivial) cofibrations}\label{gtc}

In this section we present two sets of morphisms in  $\cat{\C V}{}$, $I'$ and $J'$, which a fortiori will be sets of generating (trivial) cofibrations. Moreover, we check that almost all of the hypotheses of the recognition theorem for cofibrantly generated model categories are satisfied, see \cite[Theorem 2.1.19]{hmc}. There is however one of them which is very difficult to check, since push-outs in $\cat{\C V}{}$ are rather complicated. This last hypothesis will be verified in Section \ref{laprueba} below, after a careful scrutiny of certain push-outs.

Let us recall the following terminology from \cite{hmc}.

\begin{defn}
Let  $C$ be  a class of morphisms in a cocomplete category $\C C$:
\begin{enumerate}
\item $\inj{C}$ is the class of morphisms satisfying the right lifting property with respect to $C$, which are called \emph{$C$-injective morphisms}. 
\item $\cof{C}$ is the class of morphisms satisfying the left lifting property with respect to $\inj{C}$, called \emph{$C$-cofibrations}. 
\item $\cell{C}$ is the class of  \emph{relative $C$-cell complexes}, i.e.~transfinite compositions of push-outs of morphisms in $C$. \item If $f\colon X\r Y$ is a retract of a relative $C$-cell complex in the comma category $X\downarrow \C C$, we write $f\in\cellr{C}$. 
 \end{enumerate}
\end{defn}

\begin{rem}
The classes $\cell{C}$ and $\cof{C}$ are closed under push-outs and transfinite compositions, see \cite[Lemma 2.1.12]{hmc}, but in general $\cellr{C}$ is only closed under push-outs. The inclusions $\cell{C}\subset\cellr{C}\subset\cof{C}$ always hold. If $C$ is a set and $\C C$ is  locally presentable, then the equality $\cellr{C}=\cof{C}$ follows from the small object argument, compare \cite[Corollary 2.1.15]{hmc}. 

The classes of cofibrations, fibrations, trivial cofibrations, and trivial fibrations in $\C V$ are
$\cellr{I}=\cof{I}$, $\inj{J}$, $\cellr{J}=\cof{J}$, and $\inj{I}$, respectively.
\end{rem}

\begin{defn}
A $\C V$-functor $\varphi\colon\CC H\r \CC K$ is a \emph{local (trivial) (co)fibration} if the morphism 
$\varphi(x,y)\colon\CC H(x,y)\r \CC K(\varphi(x),\varphi(y))$ is a (trivial) (co)fibration in $\C V$ for all $x,y\in\ob\CC H$.
\end{defn}

Consider the following adjoint pair, obtained from \eqref{frees}  and \eqref{dosob}, 
$$\xymatrix@C=40pt{\C V\ar@<.5ex>[r]^-{T_{\{0,1\}}(-)_{01}}&\cat{\C V}{\{0,1\}}.\ar@<.5ex>[l]^-{\text{forget}}}$$

The first results in this section  are immediate consequences of elementary facts on adjoint functors and Grothendieck bifibrations, hence we do not include formal proofs.

\begin{prop}\label{01}\label{010}\label{011}
The following statements hold:
\begin{enumerate}
\item An $T_{\{0,1\}}(I_{01})$-injective morphism in $\cat{\C V}{}$ is the same a local trivial fibration. 

\item A $T_{\{0,1\}}(J_{01})$-injective morphism in $\cat{\C V}{}$ is the same a local fibration.

\item  A relative $T_{\{0,1\}}(I_{01})$-cell complex in $\cat{\C V}{}$ is the same a relative $T_S(I_S)$-cell complex in some fiber $\cat{\C V}{S}$. 

\item A relative $T_{\{0,1\}}(J_{01})$-cell complex in $\cat{\C V}{}$ is the same a relative $T_S(J_S)$-cell complex in some fiber $\cat{\C V}{S}$.

\item An $T_{\{0,1\}}(I_{01})$-cofibration in $\cat{\C V}{}$ is the same as a cofibration in some fiber $\cat{\C V}{S}$. 

\item  A $T_{\{0,1\}}(J_{01})$-cofibration in $\cat{\C V}{}$ is the same as a trivial cofibration in some fiber $\cat{\C V}{S}$.
\end{enumerate}
\end{prop}

Monoids in  $\C V$, such as the initial monoid $\unit$, can be regarded as $\C V$-categories with only one object. Actually, $\cat{\C V}{\{0\}}$ is isomorphic to the category $\mon{\C V}$ of monoids in $\C V$. The isomorphism is simply 
\begin{align}
\label{isomon} \cat{\C V}{\{0\}}&\st{\cong}\To\mon{\C V},\\
\nonumber \CC H&\;\mapsto\;\CC H(0,0).
\end{align}


In what follows, the symbol $\varnothing$, mostly used in this paper to denote the initial object of $\C V$, will also denote  the $\C V$-category with empty set of objects, which is initial in $\cat{\C V}{}$. There will be no possible confusion.

\begin{prop}\label{nada1}
An $\{\varnothing\r\unit\}$-injective morphism in $\cat{\C V}{}$  is the same as a $\C V$-functor surjective on objects. Moreover, a relative $\{\varnothing\r\unit\}$-cell complex is the same as an  $\{\varnothing\r\unit\}$-cofibration and the same as a cocartesian morphism $\CC K\r f_{*}\CC K$, as in \eqref{cocartesian}, associated to an injective map of sets $f$.
\end{prop}

Let 
\begin{equation}\label{iprima}
I'=T_{\{0,1\}}(I_{01})\coprod\{\varnothing\r\unit\}. 
\end{equation}

\begin{cor}\label{007}
The $I'$-injective morphisms are the local trivial fibrations surjective on objects. In particular, they are DK-equivalences.
\end{cor}

\begin{cor}\label{008}
An  $I'$-cofibration is a $\C V$-functor which factors as $\CC H\r f_{*}\CC H\r \CC K$, where the first arrow is the cocartesian morphism \eqref{cocartesian} associated to an injective map of sets $f\colon S\hookrightarrow S'$ and the second arrow is a cofibration in $\cat{\C V}{S'}$. 
\end{cor}

\begin{defn}\label{vinterval}
A \emph{$\C V$-interval} is a $\C V$-category $\CC I$ with object set $\{0,1\}$ such that $0\cong 1$ in $\pi_{0}\CC I$. A \emph{generating set of $\C V$-intervals} is a set $G$ of $\C V$-intervals such that, given a $\C V$-category $\CC K$ and two different objects $x,y\in\ob\CC K$ isomorphic $x\cong y$ in $\pi_{0}\CC K$, there exists $\CC I\in G$ and a $\C V$-functor $\phi\colon\CC I\r \CC K$ with $\phi(0)=x$ and $\phi(1)=y$.
\end{defn}

This notion of $\C V$-interval is more general than Berger--Moerdijk's \cite[Definition 1.11]{htec}. They ask $\C V$-intervals to be in addition cofibrant and, most importantly, weakly equivalent to $(\{0,1\}\onto\{0\})^*\unit$.

The following lemma is the first place were we use that $\C V$ is locally presentable.

\begin{prop}\label{lohay}
\end{prop}

\begin{proof}
There is a high enough regular cardinal $\lambda$ such that $\C V$ is locally $\lambda$-presentable, $\unit$ is $\lambda$-presentable,  and hence the functor
$$\C V\To\ho\C V\st{\pi_0}\To\operatorname{Set}$$
preserves $\lambda$-filtered colimits, see \cite[Remark 3.4 (1) and Theorem 4.1]{gbrhc}. The category $\cat{\C V}{}$ is locally $\lambda$-presentable in these circumstances, see \cite[Theorem 4.5]{vcatlplb}.  We are going to show that the set $G$ of $\C V$-intervals which are fully contained in a set of $\lambda$-presentable generators of $\cat{\C V}{}$ is generating. 

Let $\CC K$ be a $\C V$-category and $x,y\in\ob\CC K$ two different objects such that $x\cong y$ in $\pi_0\CC K$. Express $\CC K$ as an $\lambda$-filtered colimit $\CC K=\colim_{\alpha\in A}\CC K_\alpha$ of $\lambda$-presentable generators. Then $\pi_0\CC K=\colim_{\alpha\in A}\pi_0\CC K_\alpha$ in $\cat{\operatorname{Set}}{}$. In particular, there must be some $\alpha'\in A$ and objects $x',y'\in\ob\CC K_{\alpha'}$ such that $x'\cong y'$ in $\pi_0\CC K_{\alpha'}$ and the canonical functor $p\colon\CC K_{\alpha'}\r\CC K$ satisfies $p(x')=x$ and $p(y')=y$. Therefore, the full sub-$\C V$-category $\CC I$ of $\CC K_{\alpha'}$ spanned by $\{x',y'\}$ is a $\C V$-interval in $G$ (after formally renaming $x'=0$ and $y'=1$) and the restriction $\phi=p_{|_{\CC I}}\colon \CC I\r\CC K$ satisfies the requirement.
\end{proof}

We are now going to construct a set of morphisms out of a fixed set of $\C V$-intervals $G$. We will assume that $G$ is generating when necessary.

Let $i\colon\{0\}\hookrightarrow\{0,1\}$ be the inclusion. For each $\CC I\in G$, we choose a cofibrant replacement of $i^{*}\CC I$ in $\cat{\C V}{\{0\}}$,
$$\unit\into \widetilde{i^{*}\CC I}\st{\sim}\onto i^{*}\CC I.$$
Then we factor the composite $$\xymatrix{i_{*} \widetilde{i^{*}\CC I}\ar[r]& i_{*}i^{*}\CC I\ar[r]^-{\text{counit}}&\CC I}$$ into 
a cofibration 
 followed by a trivial fibration in $\cat{\C V}{\{0,1\}}$,
\begin{equation}\label{factorizacion}
i_{*}\widetilde{i^{*}\CC I}\into\tilde {\CC I}\st{\sim}\onto \CC I.
\end{equation}
We denote by
$$\theta_{\CC I}\colon \widetilde{i^{*}\CC I}\To\tilde{\CC I} $$
the composition of the cartesian mophism $\widetilde{i^{*}\CC I}\r i_{*}\widetilde{i^{*}\CC I}$ with the cofibration $i_{*}\widetilde{i^{*}\CC I}\into\tilde {\CC I}$ in \eqref{factorizacion}, and we define the set 
\begin{equation}\label{jprima}
J'=T_{\{0,1\}}(J_{01})\coprod\{\theta_{\CC I}\,;\,\CC I\in G\}. 
\end{equation}

\begin{prop}\label{thetai}
For any $\CC I\in G$, the $\C V$-functor $\theta_{\CC I}$ is a DK-equivalence and a  $I'$-cofibration.
\end{prop}

\begin{proof}
It is an $I'$-cofibration by Corollary \ref{008}. It is homotopically essentially surjective because $\CC I$ is a $\C V$-interval. Moreover, the composite
$$\xymatrix@C=15pt{\widetilde{i^{*}\CC I}(0,0)\ar[rr]^-{\theta_{\CC I}(0,0)}&&\tilde{\CC I}(0,0)\ar@{->>}[r]^{\sim}& \CC I(0,0)}$$
coincides with the trivial fibration $\widetilde{i^{*}\CC I}\st{\sim}\onto i^{*}\CC I$ under the isomorphism \eqref{isomon},  hence $\theta_{\CC I}$ is homotopically fully faithful by the 2-out-of-3 property of weak equivalences in $\C V$. Therefore, $\theta_{\CC I}$ is also a DK-equivalence.
\end{proof}

\begin{cor}
Relative $J'$-cell complexes are $I'$-cofibrations.
\end{cor}

\begin{proof}
Apart from Proposition \ref{thetai}, it is enough to notice that morphisms in $T_{\{0,1\}}(J_{01})$ are $T_{\{0,1\}}(I_{01})$-cofibrations by Proposition \ref{011}.
\end{proof}

\begin{prop}\label{unno}
In general, $I'$-injective morphisms are $J'$-injective DK-equiva\-lences. If $G$ is a generating set of $\C V$-intervals, then the converse also holds.
\end{prop}

\begin{proof}
Let us start with the first part of the statement. By Corollary \ref{007}, $I'$-injective morphisms are DK-equivalences.  By Proposition \ref{01}, $I'$-injective morphisms are $T_{\{0,1\}}(J_{01})$-injective, hence we only have to check that any $I'$-injective $\varphi\colon\CC H\r \CC K$ satisfies the left lifting property with respect to any $\theta_{\CC I}$. Consider therefore a commutative square,
$$\xymatrix{
\widetilde{i^{*}\CC I}\ar[d]_{\theta_{\CC I}}\ar[r]^-{\psi}&\CC H\ar[d]^{\varphi}\\
\tilde{\CC I}\ar[r]_-{\phi}&\CC K
}$$
Since $\varphi$ is surjective on objects, we can take $x\in \ob\CC H$ such that $\varphi(x)=\phi(1)$. Hence we can extend the previous square along the cocartesian morphism $\widetilde{i^{*}\CC I}\r i_{*}\widetilde{i^{*}\CC I}$ to
$$\xymatrix{
i_{*}\widetilde{i^{*}\CC I}\ar@{>->}[d]\ar[r]^-{\tilde \psi}&\CC H\ar[d]^{\varphi}\\
\tilde{\CC I}\ar[r]_-{\phi}&\CC K
}$$
by setting $\tilde \psi(1)=x$. Here, the left  vertical arrow is an $T_{\{0,1\}}(I_{01})$-cofibration by Proposition \ref{011}, hence a lifting $\tilde {\CC I}\r \CC H$ of the latter square exists, and it is also a lifting of the former square.

Let us prove the second part. By Proposition \ref{01}, any $J'$-injective morphism is a local fibration. Hence $J'$-injective DK-equivalences are  local trivial fibrations. By Corollary \ref{007} we just have to show that any $J'$-injective DK-equivalence $\varphi\colon\CC H\r \CC K$ is surjective on objects. Given $y\in\ob \CC K$, since $\varphi$ is a DK-equivalence, there exists $x\in\ob\CC H$ such that $\varphi(x)\cong y$ in $\pi_{0}\CC K$. Since $G$ is generating, we can choose a $\C V$-functor $\phi\colon\CC I\r\CC K$ from a $\C V$-interval $\CC I\in G$ with $\phi(0)=\varphi(x)$ and $\phi(1)=y$. We apply the lifting property in $\mon{\C V}$ to the following commutative square
$$\xymatrix{
\unit\ar@{>->}[d]\ar[rrr]^-{\id{x}}&&&\CC H(x,x)\ar@{->>}[d]^{\varphi(x,x)}_{\sim}\\
\widetilde{i^{*}\CC I}(0,0)\ar[r]_-{\theta_{\CC I}(0,0)}^-\sim&\tilde{\CC I}(0,0)\ar@{->>}[r]^-{\sim}&\CC I(0,0)\ar[r]_-{\phi(0,0)}&\CC K(\varphi(x),\varphi(x))
}$$
obtaning a $\C V$-functor $\psi\colon \widetilde{i^{*}\CC I}\To\CC H$ with $\psi(0)=x$. 
Hence, the following square commutes
$$\xymatrix{
\widetilde{i^{*}\CC I}\ar[d]_{\theta_{\CC I}}\ar[rr]^-{\psi}&&\CC H\ar[d]^{\varphi}\\
\tilde{\CC I}\ar@{->>}[r]_-{\sim}&\CC I\ar[r]_{\phi}&\CC K
}$$
Moreover, it has a lifting since $\varphi$ is $J'$-injective, hence $y$ comes from $\ob \CC H$.
\end{proof}

\section{Push-outs of free enriched functors}\label{pushpush}

This technical section, which is purely categorical, is a fundamental step towards the proof of our main theorem. Here we analize push-outs in $\cat{\C V}{S}$ of the following form
\begin{equation}\label{pussy}
\xymatrix{T_S(U_{ab})\ar[r]^-{T_S(f_{ab})}\ar[d]_{g}\ar@{}[rd]|{\text{push}}&T_S(V_{ab})\ar[d]^{g'}\\
\CC H\ar[r]_-{\varphi}& \CC K}
\end{equation}
where $a,b\in S$, see \eqref{frees} and \eqref{dosob}.

Given a (co)complete biclosed monoidal category $\C C$ with tensor product $\otimes$ and tensor unit $\unit$, we can consider \emph{monoids} in $\C C$ as well as \emph{left and right modules} over a given monoid $A$. Moreover, if $M$ and $N$ are a right and a left $A$-module with structure morphisms $\rho\colon M\otimes A\r M$ and $\lambda\colon A\otimes N\r N$, their \emph{tensor product} over $A$ is defined in the usual way, as the following coequalizer,
$$\xymatrix{M\otimes A\otimes N\ar@<-.5ex>[r]_-{\id{M}\otimes\lambda} 
\ar@<.5ex>[r]^-{\rho\otimes\id{N}} 
&
M\otimes N\ar@{->>}[r]& M\otimes_{A}N}.$$
We can also consider \emph{bimodules}, etc. The category of \emph{left $A$-modules} in $\C C$ will be denoted by $\lmod{A}$. The underlying category $\C C$ will be understood from the context. The category of \emph{right $A$-modules} is $\lmod{A^{\op}}$, where $A^{\op}$ denotes the \emph{opposite monoid}. Similarly, the category of \emph{$A$-bimodules} is $\lmod{A^{\env}}$, where $A^{\env}=A\otimes A^{\op}$ is the \emph{enveloping monoid}.

The forgetful functor $\lmod{A}\r\C C$ preserves colimits since it has a right adjoint $\hom_{\C C}^l(A,-)$. Here, given an object $X$ in $\C C$, $\hom_{\C C}^l(X,-)$ denotes the right adjoint of $X\otimes-$. Therefore colimits in module categories are computed in $\C C$.

If $\CC H$ is a $\C V$-category, endomorphism objects $\CC H(x,x)$ are monoids in $\C V$,  the morphism object $\CC H(x,y)$ is a left-$\CC H(y,y)$-right-$\CC H(x,x)$-bimodule, and composition factors through the tensor product over the endomorphism monoid of the middle object
$$\xymatrix{
\CC H(y,z)\otimes\CC H(x,y)\ar[rr]^-{c_{\CC H}(x,y,z)}\ar@{->>}[rd]&& \CC H(x,z)\\
&\CC H(y,z)\otimes_{\CC H(y,y)}\CC H(x,y)\ar[ru]_-{\quad\bar c_{\CC H}(x,y,z)}&}$$
All these arrows are  left-$\CC H(z,z)$-right-$\CC H(x,x)$-bimodule morphisms. Moreover, the morphism $\bar c_{\CC H}(x,y,z)$, called \emph{reduced composition morphism (interpolating at $y$)}, is an isomorphism if either $x=y$ or $y=z$.

Denote by $\mon{\C C}$ the category of monoids in $\C C$. The forgetful functor has a left adjoint, the \emph{free monoid functor}~$T$,
$$\xymatrix{\C C\ar@<.5ex>[r]^-{T}&\mon{\C C}\ar@<.5ex>[l]^-{\text{forget}}}$$ 
defined by the usual tensor algebra construction
$$T(U)=\coprod_{n\geq 0} U^{\otimes n}.$$
Multiplication is given by the distributibity of $\coprod$ over $\otimes$, and the monoid unit of $T(U)$ is the inclusion of the factor $n=0$ of the coproduct, which is $U^{\otimes 0}=\unit$. 
The unit of the adjunction $U\r T(U)$ is the inclusion of the factor $n=1$ of the coproduct, and the counit $T(A)\r A$ is defined by the multiplication and the unit of the monoid $A$. 

The \emph{category $\mor{\C{C}}$ of morphisms in $\C C$} can be regarded as the category of functors $\dos\r\C{C}$, where $\dos$ is the category with two objects, $0$ and $1$, and only one non-identity morphism $0\r 1$, i.e.~it is the poset $\{0<1\}$. A morphism $f\colon U\r V$ in $\C C$ is identified with the functor $f\colon\dos\r\C{C}$ defined by $f(0)=U$, $f(1)=V$ and $f(0\r 1)=f$.

The category $\mor{\C{C}}$ carries a biclosed monoidal structure given by the \emph{push-out product} of morphisms $f\odot g$,
$$\xy/r2.3pt/:
(0,0)*+{U\otimes X}="a",
(35,0)*+{V\otimes X}="b",
(0,-18.2)*+{U\otimes Y}="c",
(35,-20)*+{U\otimes Y\!\bigcup\limits_{U\otimes X}\! V\otimes X}="d",
(70,-35)*+{V\otimes Y}="e"
\ar"a";"b"^{f\otimes \id{X}}
\ar@{}"a";"d"|-{\text{push}}\ar"a";"c"_{\id{U}\otimes g}
\ar"b";"d"
\ar"c";(17,-18.2)
\ar@/^15pt/"b";"e"^{\id{V}\otimes g}
\ar@/_15pt/"c";"e"_{f\otimes\id{Y}}
\ar(40,-23);"e"^-{f\odot g}
\endxy$$
Notice that $f\odot g$ is an isomorphism if $f$ or $g$ is. 
This monoidal structure is symmetric provided  $\otimes$ is. If $\varnothing$ denotes the initial object of $\C C$, the functor 
\begin{align*}
\C{C}&\To\mor{\C C},\\
X&\;\mapsto\; (\varnothing\r X),
\end{align*}
is strong (symmetric) monoidal. Moreover, $f\odot(\varnothing\r X)=f\otimes X$ and $(\varnothing\r X)\odot f=X\otimes f$. 


\begin{lem}[{\cite[Lemma 4.1]{htnso}}]\label{podot}
Given two push-out diagrams in $\C{C}$, $i=1,2$,
$$\xymatrix{U_i\ar[r]^{f_i}\ar[d]_{g_i}\ar@{}[rd]|{\text{push}}&V_i\ar[d]^{g'_i}\\X_i\ar[r]_{f'_i}&Y_i}$$
the following diagram in $\C{C}$ is also a push-out,
$$\xy
(0,-2)*{U_1\otimes V_2\!\!\!\!\bigcup\limits_{U_1\otimes U_2}\!\!\!\!V_1\otimes U_2},
(0,-22)*{X_1\otimes Y_2\!\!\!\!\bigcup\limits_{X_1\otimes X_2}\!\!\!\!Y_1\otimes X_2},
(40,0)*{V_1\otimes V_2},
(40,-20)*{Y_1\otimes Y_2},
(20,-10)*{\text{\scriptsize push}}
\ar(16,0);(33,0)^-{f_1\odot f_2}
\ar(40,-3);(40,-17)^-{g_1'\otimes g_2'}
\ar(16,-20);(33,-20)_-{f'_1\odot f'_2}
\ar(0,-6);(00,-17)_-{g_1\otimes g_2' \!\!\!\!\bigcup\limits_{g_1\otimes g_2}\!\!\!\! g_1'\otimes g_2}
\endxy$$
\end{lem}

Given morphisms $f_{i}\colon U_{i}\r V_{i}$ in $\C C$, $1\leq i\leq n$, the target of $f_{1}\odot\cdots\odot f_{n}$ is the iterated tensor product of the targets $V_{1}\otimes\cdots\otimes V_{n}$. This object is the colimit of the diagram
$$f_{1}\otimes\cdots\otimes f_{n}\colon \dos^{n}\To\C{C},$$
since $\dos^{n}$ has a final object $(1,\st{n}\dots,1)$. The source of $f_{1}\odot\cdots\odot f_{n}$ is the colimit of the restriction of this diagram to the full subcategory of $\dos^{n}$ obtained by removing the final object, and the push-out product
$$f_{1}\odot\cdots\odot f_{n}\colon \colim_{\dos^n\setminus\{(1,\st{n}\dots,1)\}}f_{1}\otimes\cdots\otimes f_{n}\To \colim_{\dos^{n}}f_{1}\otimes\cdots\otimes f_{n}=V_{1}\otimes\cdots\otimes V_{n}$$
is the morphism induced by the inclusion $\dos^n\setminus\{(1,\st{n}\dots,1)\}\subset\dos^n$. 
Later, in diagrams, we will often drop the source and target of such a push-out product from notation, and simply write
$$\xymatrix@C=45pt{\bullet\ar[r]^-{f_{1}\odot\cdots\odot f_{n}}&\bullet.}$$

The universal property of the source $\colim_{\dos^n\setminus\{(1,\st{n}\dots,1)\}}f_{1}\otimes\cdots\otimes f_{n}$  refers to canonical morphisms in $\C C$ 
$$\kappa_{i}\colon V_{1}\otimes\cdots\otimes V_{i-1}\otimes U_{i}\otimes V_{i+1}\otimes\cdots\otimes V_{n}\To
\colim_{\dos^n\setminus\{(1,\st{n}\dots,1)\}}f_{1}\otimes\cdots\otimes f_{n},\quad 1\leq i\leq n,$$
with $(f_{1}\odot\cdots\odot f_{n})\kappa_{i}=\id{}^{\otimes(i-1)}\otimes f_{i}\otimes\id{}^{\otimes(n-i)}$.
Any collection of morphisms from the objects obtained from  $V_{1}\otimes\cdots\otimes V_{n}$ by replacing the target of each push-out product factor with the source 
$$g_{i}\colon V_{1}\otimes\cdots\otimes V_{i-1}\otimes U_{i}\otimes V_{i+1}\otimes\cdots\otimes V_{n}\To
X,\quad 1\leq i\leq n,$$
such that the following squares commute, $1\leq i<j\leq n$,
\begin{equation}\label{beber}
\xymatrix@C=50pt{
V_{1}\otimes\cdots\otimes U_{i}\otimes \cdots\otimes U_{j}\otimes \cdots\otimes V_{n}
\ar[r]^-{
\begin{array}{c}
\scriptstyle
\id{}\otimes\cdots\otimes f_{i}\otimes\cdots\otimes\id{}
\\
\vspace{-10pt}
\end{array}
}
\ar[d]_-{
\id{}\otimes\cdots\otimes f_{j}\otimes\cdots\otimes\id{}
}&
V_{1}\otimes\cdots\otimes V_{i}\otimes \cdots\otimes U_{j}\otimes \cdots\otimes V_{n}\ar[d]^-{g_{j}}\\
V_{1}\otimes\cdots\otimes U_{i}\otimes \cdots\otimes V_{j}\otimes \cdots\otimes V_{n}\ar[r]_-{g_{i}}&
X
}
\end{equation}
induces a unique morphism $g\colon \colim_{\dos^n\setminus\{(1,\st{n}\dots,1)\}}f_{1}\otimes\cdots\otimes f_{n} \r X$ such that $g_{i}=g\kappa_{i}$, $1\leq i\leq n$. Compare the paragraph preceding \cite[Lemma 2.2]{hcctaimc}.

The push-out of a diagram of monoids in $\C C$
$$\xymatrix{T(U)\ar[r]^-{T(f)}\ar[d]_{g}\ar@{}[rd]|{\text{push}}&T(V)\ar[d]^{g'}\\
A\ar[r]_-{\varphi}& B}$$
can be constructed as follows \cite[proof of Lemma
6.2]{ammmc}. Let $\bar g\colon U\r M$ be the adjoint of $g$. The underlying morphism of $\varphi$ is a transfinite (countable) composition in $\C C$
$$A=B_{0}\st{\varphi_{1}}\To B_{1}\r\cdots\r B_{t-1}\st{\varphi_{t}}\To B_{t}\r\cdots$$
such that the \emph{bonding morphism} $\varphi_{t}$ is a push-out in $\C C$ of the form
$$\xymatrix@C=40pt{\bullet\ar[r]^-{(A\otimes f)^{\odot t}\otimes A}\ar[d]_{\psi_{t}}\ar@{}[rd]|{\text{push}}&\bullet\ar[d]^{\bar\psi_{t}}\\
B_{t-1}\ar[r]_-{\varphi_{t}}&B_{t}}$$
Notice that the target of the upper horizontal arrow is
\begin{equation}\label{target}
A\otimes V\otimes \cdots\otimes A \otimes V\otimes A,
\end{equation}
with $t$ copies of $V$ between $t+1$ copies of $A$. The morphism $\psi_{t}$ is defined by the $t$ maps
$$\xymatrix{A\otimes  \cdots\otimes A\otimes U\otimes A\otimes\cdots\otimes A\ar[d]_{\id{}\otimes \bar g\otimes\id{}}\\
A\otimes  \cdots\otimes \underbrace{A\otimes A\otimes A}\otimes\cdots\otimes A\ar[d]_{\id{}\otimes\text{multiplication}\otimes\id{}}\\
A\otimes  \cdots\otimes A\otimes\cdots\otimes A\ar[d]_{\bar\psi_{t-1}}^{\text{or the identity if }t=1}\\
B_{t-1}
}$$
where $U$ replaces one of the   copies of $V$ in \eqref{target}. The multiplication
$B\otimes B\r B$ is the colimit of the unique morphisms
$$c_{s,t}\colon B_{s}\otimes B_{t}\To B_{s+t}$$
such that
\begin{align*}
c_{s,t}(\varphi_{s}\otimes\id{})&=\varphi_{s+t}c_{s-1,t},&
c_{s,t}(\id{}\otimes\varphi_{t})&=\varphi_{s+t}c_{s,t-1},
\end{align*}
and $c_{s,t}(\bar\psi_{s}\otimes\bar\psi_{t})$ is the composite
$$\xymatrix{
\hspace{-11pt}(A\otimes V)^{\otimes s}\otimes \underbrace{A\otimes (A}\otimes V)^{\otimes t}\otimes A\ar[d]_{\id{}\otimes\text{multiplication}\otimes\id{}}\\
\hspace{-16pt}(A\otimes V)^{\otimes s}\otimes (A\otimes V)^{\otimes t}\otimes A
\ar[d]_{\bar\psi_{s+t}}\\
B_{s+t}
}$$
This uses the push-out definition of $B_{s}\otimes B_{t}$ given by Lemma \ref{podot}. 
The unit  of $B$  is the composite $$\unit\st{\text{unit}}\To A\st{\varphi}\To B$$ and the adjoint of $g'$ is  
$$\xymatrix{V\cong \unit\otimes V\otimes \unit\ar[rr]^-{\text{unit}\otimes\id{V}\otimes\text{unit}}&& A\otimes V\otimes A\ar[r]^-{\bar\psi_{1}}& B_{1}\ar[r]& B.}$$

Taking $\C C=\graph{\C V}{S}$, we can describe the push-out \eqref{pussy} in $\cat{\C V}{S}$
 in the following way. Given $x,y\in S$, $\varphi(x,y)\colon \CC H(x,y)\r  \CC K(x,y)$ is a transfinite composition in $\C V$
$$\xymatrix{\CC H(x,y)=\CC K(x,y)_{0}\ar[r]^-{\varphi(x,y)_{1}} & \CC K(x,y)_{1}\r\cdots\r \CC K(x,y)_{t-1}\ar[r]^-{\varphi(x,y)_{t}} & \CC K(x,y)_{t}\r\cdots}$$
such that the bonding morphism $\varphi(x,y)_{t}$ is a push-out in $\C V$ of the form
\begin{equation}\label{cosaque}
\xymatrix@C=100pt{\bullet\ar[r]^-{\CC H(b,y)\otimes( f\otimes \CC H(b,a))^{\odot (t-1)}\odot f\otimes \CC H(x,a)}\ar[d]_{\psi(x,y)_{t}}\ar@{}[rd]|{\text{push}}&\bullet\ar[d]^{\bar\psi(x,y)_{t}}\\
\CC K(x,y)_{t-1}\ar[r]_-{\varphi(x,y)_{t}}&\CC K(x,y)_{t}}
\end{equation}
The target of the upper horizontal arrow is
\begin{equation}\label{target2}
\CC H(b,y)
\otimes
V\otimes\CC H(b,a)\otimes \cdots
\otimes V\otimes\CC H(b,a)\otimes
V
\otimes 
\CC H(x,a),
\end{equation}
with $t$ copies of $V$, 
and $\psi(x,y)_{t}$ is defined by the maps
\begin{equation}\label{pegando}
\xymatrix{\CC H(b,y)
\otimes
V\otimes\cdots
\otimes
\CC H(b,a)
\otimes
U\otimes\CC H(b,a)
\otimes
\cdots
\otimes
V
\otimes 
\CC H(x,a)\ar[d]_{\id{}\otimes \bar g\otimes\id{}}\\
\CC H(b,y)
\otimes
V\otimes\cdots
\otimes
\underbrace{\CC H(b,a)
\otimes
\CC H(a,b)\otimes\CC H(b,a)}
\otimes
\cdots
\otimes
V
\otimes 
\CC H(x,a)\ar[d]_{\id{}\otimes\text{composition}\otimes\id{}}\\
\CC H(b,y)
\otimes
V\otimes\cdots
\otimes
\CC H(b,a)
\otimes
\cdots
\otimes
V
\otimes 
\CC H(x,a)\ar[d]_{\bar\psi(x,y)_{t-1}}^{\text{or the identity if }t=1}\\
\CC K(x,y)_{t-1}}
\end{equation}
Here $\bar g\colon U\r \CC H(a,b)$ is the adjoint of $g$ and $U$ replaces one of the copies of $V$ in the middle. If it replaces the first or last occurrence of $V$, the map is similarly defined. This decomposition of $\varphi(x,y)\colon \CC H(x,y)\r  \CC K(x,y)$ in $\C V$ is actually a decomposition in the category of left-$\CC H(y,y)$-right-$\CC H(x,x)$-bimodules since all arrows are clearly bimodule morphisms.

Composition
$$c_{\CC K}(x,y,z)\colon\CC K(y,z)\otimes \CC K(x,y)\To \CC K(x,z)$$ 
is defined as the colimit of maps
$$c_{\CC K}(x,y,z)_{s,t}\colon\CC  K(y,z)_{s}\otimes\CC K(x,y)_{t}\To\CC  K(x,z)_{s+t}$$ 
such that $c_{\CC K}(x,y,z)_{0,0}=c_{\CC H}(x,y,z)$ is composition in $\CC H$,
\begin{align*}
c_{\CC K}(x,y,z)_{s,t}(\varphi(y,z)_{s}\otimes\id{})&=\varphi(x,z)_{s+t}c_{\CC K}(x,y,z)_{s-1,t},\\
c_{\CC K}(x,y,z)_{s,t}(\id{}\otimes\varphi(x,y)_{s})&=\varphi(x,z)_{s+t}c_{\CC K}(x,y,z)_{s,t-1},
\end{align*}
and $c(x,y,z)_{s,t}(\bar\psi(y,z)_{s},\bar\psi(x,y)_{t})$ is
$$\xy
(0,0)*+{\CC H(b,z)
\hspace{-2pt}\otimes\hspace{-2pt}
(V\hspace{-2pt}\otimes\hspace{-2pt}\CC H(b,a))^{\otimes(s-1)}
\hspace{-2pt}\otimes\hspace{-2pt}
V
\hspace{-2pt}\otimes\hspace{-2pt} 
\underbrace{\CC H(y,a)
\hspace{-2pt}\otimes\hspace{-2pt}
\CC H(b,y)}
\hspace{-2pt}\otimes\hspace{-0pt}
(V\hspace{-2pt}\otimes\hspace{-2pt}\CC H(b,a))^{\otimes(t-1)}
\hspace{-2pt}\otimes\hspace{-2pt}
V
\hspace{-2pt}\otimes\hspace{-2pt} 
\CC H(x,a)}="a"
,
(0,-15)*+{\CC H(b,z)
\hspace{-2pt}\otimes\hspace{-2pt}
(V\hspace{-2pt}\otimes\hspace{-2pt}\CC H(b,a))^{\otimes(s-1)}
\hspace{-2pt}\otimes\hspace{-2pt}
V
\hspace{-2pt}\otimes\hspace{-2pt} 
\CC H(b,a)
\hspace{-2pt}\otimes\hspace{-2pt}
(V\hspace{-2pt}\otimes\hspace{-2pt}\CC H(b,a))^{\otimes(t-1)}
\hspace{-2pt}\otimes\hspace{-2pt}
V
\hspace{-2pt}\otimes\hspace{-2pt} 
\CC H(x,a)}="b",
(0,-30)*+{\CC K(x,z)_{s+t}}="c"
\ar"a";"b"_{\id{}\otimes\text{composition}\otimes\id{}}\ar"b";"c"_{\bar\psi(x,z)_{s+t}}
\endxy
$$
Identities are given by the compositions
$$\xymatrix{\unit\ar[r]^-{\id{x}}&\CC H(x,x)\ar[r]^-{\varphi(x,x)}& \CC K(x,x)}$$
and the adjoint of $g'$ is
$$\xymatrix{V\cong \unit\otimes V\otimes \unit\ar[rr]^-{\id{b}\otimes\id{V}\otimes\id{a}}&& \CC H(b,b)\otimes V\otimes \CC H(a,a)\ar[r]^-{\bar\psi(a,b)_{1}}& \CC K(a,b)_{1}\ar[r]& \CC K(a,b).}$$
Here we denote the identity morphism $\id{V}\colon V\r V$ in the same way as the identities $\id{a}\colon\unit\r\CC H(a,a)$ and $\id{b}\colon\unit\r\CC H(b,b)$ in the $\C V$-category $\CC H$.

%
%

The previous decomposition of $\varphi(x,y)\colon\CC H(x,y)\r \CC K(x,y)$ as a transfinite composition of push-outs in $\C V$ is nice and detailed, and will later be used to deduce homotopical properties of $\varphi$ when $f$ is a cofibration, etc. However this is not enough. We need to be able to say something (homotopical) about $\varphi(x,y)$ as a left $\CC H(y,y)$-module morphism, as a right $\CC H(x,x)$-module morphism, and even as a monoid morphism when $x=y$. Moreover, we must be able to deduce homotopical properties of reduced composition morphisms in $\CC K$ out of the properties of reduced compositions in $\CC H$. This is why we need similar decompositions of $\varphi(x,y)$ in module (and monoid) categories, and also of the commutative square $\bar c_{\CC H}\r\bar c_{\CC K}$ induced by $\varphi$, see \eqref{composdiag1}. The rest of the section is devoted to this. The statements, long and detailed, present the new decompositions of $\varphi(x,y)$ (and $\bar c_{\CC H}\r\bar c_{\CC K}$), 
very similar in 
flavour to the one above. The idea is that, if the monoid $\CC H(z,z)$ appears in the statement, then morphism objects in the upper arrow of \eqref{cosaque} must be replaced with the corresponding reduced composition in $\CC H$ interpolating at $z$ (unless the reduced composition is an isomorphism). We have tried to keep proofs short and not too detailed. The proof of the third lemma is actually not included since it would be very similar to the previous ones.

The category $\lmod{A^{\env}}$ of bimodules over a monoid $A$ is a (non-symmetric) monoidal category with tensor product $\otimes_{A}$, so we can consider monoids in $\lmod{A^{\env}}$. The category $\mon{\lmod{A^{\env}}}$ is isomorphic to the comma category $A\downarrow\mon{\C V}$ of monoids under $A$. Indeed, a monoid morphism $A\r B$ induces an $A$-bimodule structure on $B$ and the multiplication $B\otimes B\r B$ factors through an $A$-bimodule morphism $B\otimes_{A} B\r B$. Conversely, the monoid $A$ is initial in $\mon{\lmod{A^{\env}}}$, so the unit $A\r B$ of any monoid $B$ in $\lmod{A^{\env}}$ is an object in $A\downarrow\mon{\C V}$.

\begin{lem}\label{monoid}
Given $x\in S$, the monoid morphism $\varphi(x,x)\colon \CC H(x,x)\r\CC K(x,x)$ in the push-out \eqref{pussy} is a transfinite composition in $\mon{\lmod{\CC H(x,x)^{\env}}}$
$$\CC H(x,x)=\CC K_{x,0}\st{\varphi'_{x,1}}\To \CC K_{x,1}\r\cdots\r \CC K_{x,t-1}\st{\varphi'_{x,t}}\To \CC K_{x,t}\r\cdots$$
such that the bonding morphism 
$\varphi'_{x,t}$ fits into a push-out in that category
\begin{equation}\label{monoiddiag}
\xymatrix@C=150pt{
\bullet
\ar[r]^-{T(
\CC H(b,x)\otimes (f\odot \bar c_{\CC H}(b,x,a))^{\odot(t-1)}\odot f\otimes\CC H(x,a)
)}
\ar[d]_{\psi'_{x,t}}\ar@{}[rd]|{\text{push}}&
\bullet
\ar[d]^{\bar\psi'_{x,t}}\\
\CC K_{x,t-1}\ar[r]_-{\varphi'_{x,t}}&\CC K_{x,t}}
\end{equation}
The target of the upper horizontal arrow is the free monoid on the $\CC H(x,x)$-bimodule \eqref{target2} for $y=x$. The morphism $\psi'_{x,t}$ is defined by a map from the $\CC H(x,x)$-bimodule which generates the free monoid in the source. This generating bimodule is the source of a push-out product, so the previous map is defined by maps from the $\CC H(x,x)$-bimodules which are obtained from \eqref{target2} by replacing the target of each push-out product factor with the source. If we replace a $V$ with a $U$, the map is defined as in \eqref{pegando}, putting $y=x$ and changing  $\bar\psi(x,y)_{t-1}$   for $\bar\psi'_{x,t-1}$. 
%
If we replace the $s^{\text{th}}$ copy of $\CC H(b,a)$ with $\CC H(x,a)\otimes_{\CC H(x,x)} \CC H(b,x)$, $1\leq s\leq t-1$, then the map is
$$\xymatrix{\CC H(b,x)
\otimes \cdots
\otimes V\otimes\CC H(x,a)\otimes_{\CC H(x,x)} \CC H(b,x)\otimes
V
\otimes \cdots\otimes
\CC H(x,a)\ar[d]_{\bar\psi'_{x,s}\otimes_{\CC H(x,x)}\bar\psi'_{x,t-s}}\\
\CC K_{x,s}\otimes_{\CC K_{x,0}} \CC K_{x,t-s}
\ar[d]_{
\text{bonding}\otimes \text{bonding}
}
\\
\CC K_{x,t-1}\otimes_{\CC K_{x,0}} \CC K_{x,t-1}\ar[d]_{\text{multiplication}}\\
\CC K_{x,t-1}\\
}$$
\end{lem}

\begin{proof}
First of all, notice that the maps in the statement satisfy the appropriate compatibility conditions in order to define a morphism from the source of the free monoid morphism on the push-out product. This follows easily by induction. Indeed, if we replace two consecutive copies of $V$  the appropriate square, like \eqref{beber}, commutes by the associativity of composition in $\CC H$ and the induction hypothesis. If the  two copies of $V$ are not consecutive we simply use the functoriality of $\otimes$ instead of associativity. If we replace two copies of $\CC H(b,a)$, the commutativity follows from the associativity of the multiplication in $\CC K_{x,t-1}$. If we replace a copy of $V$ and a copy of $\CC H(b,a)$, it also follows easily by induction. Hence, the monoids $\CC K_{x,t}$ and the monoid morphisms $\varphi_{x,t}'$ are well defined. In particular the colimit $\CC K_{x}=\colim_{t\geq 0}\CC K_{x,t}$ and the transfinite composition $\varphi_{x}'\colon \CC H(x,x)\r\CC K_{x}$ are well defined.

We now define two mutually inverse monoid isomorphisms
$$
\xymatrix{\CC K(x,x)\ar@<.5ex>[r]^-{\Phi_{x}}&\CC K_{x}\ar@<.5ex>[l]^-{\Phi'_{x}}}
$$
with $\varphi(x,x)=\Phi'_x\varphi'_x$ and $\varphi'_x=\Phi_x\varphi(x,x)$.

The morphism $\Phi'_{x}$ is defined by applying the universal property of a colimit to a sequence of compatible monoid morphisms $\Phi'_{x,t}\colon \CC K_{x,t}\r\CC K(x,x)$, $t\geq 0$, each of which is defined by applying the universal property of a push-out of monoids in $\lmod{\CC H(x,x)^{\env}}$ to the commutative square
$$\xymatrix@C=150pt{
\bullet
\ar[r]^-{T(
\CC H(b,x)\otimes (f\odot \bar c_{\CC H}(b,x,a))^{\odot(t-1)}\odot f\otimes\CC H(x,a)
)}
\ar[d]_{\psi'_{x,t}}&
\bullet
\ar[d]^{\Psi'_{x,t}}\\
\CC K_{x,t-1}\ar[r]_-{\Phi'_{x,t-1}}&\CC K(x,x)}$$
where $\Psi'_{x,t}$ is defined by $\bar\psi(x,x)_{t}$. The initial step is taken to be the monoid morphism $\Phi'_{x,0}=\varphi(x,x)\colon \CC K_{x,0}=\CC H(x,x)\r\CC K(x,x)$. The commutativity of this square also follows by induction in a straightforward way, using the definition of the multiplication in $\CC K(x,x)$ and the fact that $\Phi'_{x,0},\dots,\Phi'_{x,t-1}$ are monoid morphisms.

Similarly, $\Phi_{x}$ is defined by applying the universal property of a colimit to a sequence of compatible morphisms $\Phi_{x,t}\colon \CC K(x,x)_{t}\r\CC K_{x}$, $t\geq 0$, in $\C V$, each of which is defined by applying the universal property of a push-out in $\C V$ to the commutative square
$$\xymatrix@C=100pt{\bullet\ar[r]^-{\CC H(b,x)( f\otimes \CC H(b,a))^{\odot (t-1)}\odot f\otimes \CC H(x,a)}\ar[d]_{\psi(x,x)_{t}}\ar@{}[rd]|{\text{push}}&\bullet\ar[d]^{\Psi_{x,t}}\\
\CC K(x,x)_{t-1}\ar[r]_-{\Phi_{x,t-1}}&\CC K_{x}}$$
where $\Psi_{x,t}$ is defined by $\bar\psi'_{x,t}$.
The initial step $\Phi_{x,0}\colon \CC K(x,x)_{0}=\CC H(x,x)\r\CC K_{x}$ is taken to be $\varphi'_{x}$.  Again,  commutativity  follows  straightforwardly by induction. 

The morphism $\Phi_{x}\colon\CC K(x,x)\r\CC K_x$ in $\C V$ is actually a monoid morphism under $\CC H(x,x)$. Indeed, it easily follows by induction that the square
$$\xymatrix@C=50pt{\CC  K(x,x)_{s}\otimes\CC K(x,x)_{t}\ar[r]^-{c_{\CC K}(x,x,x)_{s,t}}\ar[d]_{\Phi_{x,s}\otimes \Phi_{x,t}}&\CC  K(x,x)_{s+t}\ar[d]^{\Phi_{x,s+t}}\\
\CC K_{x}\otimes\CC K_{x}\ar[r]^-{\text{mult.}}&\CC K_{x}}$$
commutes for all $s,t\geq 0$.

The equations $\Phi_x\Phi'_x=\id{\CC K_{x}}$ and $\Phi'_x\Phi_x=\id{\CC K(x,x)}$ now follow from the universal properties of colimits.
\end{proof}

In order to understand Lemma  \ref{monoid} better it is illustrative to look at what happens in special cases. If $a=x$ or $b=x$ then the reduced composition $\bar c_{\CC H}(b,x,a)$ is an isomorphism, so $\varphi'_{x,t}$ in \eqref{monoiddiag} is an isomorphism for $t>1$. This means that $\CC K(x,x)$ can be obtained from $\CC H(x,x)$ as a single push-out in the category of monoids in $\lmod{\CC H(x,x)^{\env}}$. If moreover $a=b=x$ then in \eqref{monoiddiag} for $t=1$, $T$ is applied to the free bimodule morphism on $f$, hence we simply have a push-out in $\mon{\C V}$
$$\xymatrix{
T(U)\ar[r]\ar[d]_{g(x,x)}\ar@{}[rd]|{\text{push}}&T(V)\ar[d]\\
\CC H(x,x)\ar[r]&\CC K(x,x)
}$$

\begin{lem}\label{composition}
Consider the push-out \eqref{pussy}. The morphism $\xi_{x,y,z}$ in the following diagram in $\C V$ 
\begin{equation}\label{composdiag1}
\xymatrix{
\CC H(y,z)\otimes_{\CC H(y,y)}\CC H(x,y)\ar[r]^-{\bar c_{\CC H}(x,y,z)}
\ar[d]_{\varphi(y,z)\otimes \varphi(x,y)}\ar@{}[rd]|{\text{push}}
&\CC H(x,z)\ar@(r,u)[rdd]^-{\varphi(x,z)}\ar[d]^{\tilde \varphi}\\
\CC K(y,z)\otimes_{\CC K(y,y)}\CC K(x,y)\ar@(d,l)[rrd]_{\bar c_{\CC K}(x,y,z)}
\ar[r]_-{\tilde c}
&\CC K(x,y,z)_{0}
\ar[rd]^{\xi_{x,y,z}}\\
&&\CC K(x,z)
}
\end{equation}
is a transfinite composition 
$$
\xymatrix@C=35pt{
\CC K(x,y,z)_{0}\ar[r]^-{\xi_{x,y,z,1}}& \CC K(x,y,z)_{1}\r\cdots\r \CC K(x,y,z)_{t-1}\ar[r]^-{\xi_{x,y,z,t}}& \CC K(x,y,z)_{t}\r\cdots}
$$
such that the bonding morphism 
$\xi_{x,y,z,t}$ fits into a push-out in $\C V$ 
\begin{equation}\label{composdiag2}
\xymatrix@C=150pt{
\bullet
\ar[r]^-{
\bar c_{\CC H}(b,y,z)\odot (f\odot \bar c_{\CC H}(b,y,a))^{\odot(t-1)}\odot f\odot\bar c_{\CC H}(x,y,a)
}
\ar[d]_{\psi(x,y,z)_{t}}\ar@{}[rd]|{\text{push}}&
\bullet
\ar[d]^{\bar\psi(x,y,z)_{t}}\\
\CC K(x,y,z)_{t-1}\ar[r]_-{\xi_{x,y,z,t}}&\CC K(x,y,z)_{t}}
\end{equation}
The target of the upper horizontal arrow is \eqref{target2} for $y=z$. The morphism $\psi(x,y,z)_{t}$ 
departs from the source of a push-out product, so it is defined by maps from the objects obtained from \eqref{target2} by replacing the target of each push-out product factor with the source.  If we replace a $V$ with a $U$, the map is defined as in \eqref{pegando}, putting $y=z$, changing $\bar\psi(x,y)_{t-1}$
for $\bar\psi(x,y,z)_{t-1}$ if $t>1$, and the identity for $\tilde\varphi$ if $t=1$. If we replace the $s^{\text{th}}$ copy of $\CC H(b,a)$ with $\CC H(x,a)\otimes_{\CC H(x,x)} \CC H(b,x)$, $1\leq s\leq t-1$, then the map is
\begin{equation}\label{registroauxiliar}
\xymatrix{\CC H(b,z)
\otimes \cdots
\otimes V\otimes\CC H(y,a)\otimes_{\CC H(y,y)} \CC H(b,y)\otimes
V
\otimes \cdots\otimes
\CC H(x,a)\ar[d]_{\bar\psi(y,z)_{s}\otimes_{\CC H(y,y)}\bar\psi(x,y)_{t-s}}\\
\CC K(y,z)_{s}\otimes_{\CC H(y,y)} \CC K(x,y)_{t-s}\ar[d]^{\text{tensor product of canonical maps to the colimit}}
\\
\CC K(y,z)\otimes_{\CC H(y,y)} \CC K(x,y)\ar@{->>}[d]^{\text{quotient map to a (bigger) coequalizer}}\\
\CC K(y,z)\otimes_{\CC K(y,y)} \CC K(x,y)\ar[d]_{\tilde c}\\
\CC K(x,y,z)_{0}\ar[d]_{
\text{bonding}
}^-{\text{or the identity if }t=1}\\
\CC K(x,y,z)_{t-1}
}
\end{equation}
If we replace $\CC H(b,z)$ with $\CC H(y,z)\otimes_{\CC H(y,y)}\CC H(b,y)$ or $\CC H(x,a)$ with $\CC  H(y,a)\otimes_{\CC H(y,y)}\CC H(x,y)$, then the maps are the composion of 
$$\xymatrix@C=63pt{\CC H(y,z)\otimes_{\CC H(y,y)}\CC H(b,y)\otimes V\otimes\cdots \otimes \CC H(x,a)\ar[r]^-{\id{}\otimes_{\CC H(y,y)}\bar\psi(x,y)_{t}}&
\CC K(y,z)_{0}\otimes_{\CC H(y,y)} \CC K(x,y)_{t}
}$$
$$\xymatrix@C=63pt{\CC H(b,z)\otimes\cdots \otimes V\otimes \CC  H(y,a)\otimes_{\CC H(y,y)}\CC H(x,y)\ar[r]^-{\bar\psi(y,z)_{t}\otimes_{\CC H(y,y)}\id{}}&
\CC K(y,z)_{t}\otimes_{\CC H(y,y)} \CC K(x,y)_{0}
}$$
with the maps $\CC K(y,z)_{0}\otimes_{\CC H(y,y)} \CC K(x,y)_{t}\r\CC K(x,y,z)_{t-1}\leftarrow\CC K(y,z)_{t}\otimes_{\CC H(y,y)} \CC K(x,y)_{0}
$ defined as in \eqref{registroauxiliar}, respectively.
\end{lem}

\begin{proof}
As in the proof of Lemma \ref{monoid}, the morphisms in the statement satisfy the appropriate compatibility conditions in order to define a morphism from the source of the push-out product. This follows easily by induction. Let us just point out that, if we replace two copies of $\CC H(b,a)$, the appropriate square commutes by the relations in $\CC K(y,z)\otimes_{\CC K(y,y)}\CC K(x,y)$ imposed by the coequalizer definition of $\otimes_{\CC K(y,y)}$. The same happens if we replace the copy of $\CC H(x,a)$ and a copy of $\CC H(b,a)$, etc. In particular the colimit $\CC K(x,y,z)=\colim_{t\geq 0}\CC K(x,y,z)_{t}$ and the transfinite composition $\xi'_{x,y,z}\colon \CC K(x,y,z)_0\r\CC K(x,y,z)$ of the $\xi_{x,y,z,t}$ are well defined.

We now define two mutually inverse isomorphisms in $\C V$
$$
\xymatrix{\CC K(x,z)\ar@<-.5ex>@{<-}[r]_-{\Phi'_{x,y,z}}&\CC K(x,y,z)\ar@<-.5ex>@{<-}[l]_-{\Phi_{x,y,z}}}
$$
with $\xi'_{x,y,z}=\Phi_{x,y,z}\xi_{x,y,z}$ and $\xi_{x,y,z}=\Phi'_{x,y,z}\xi'_{x,y,z}$.

The morphism $\Phi'_{x,y,z}$ is defined by applying the universal property of a colimit to a sequence of compatible morphisms $\Phi_{x,y,z,t}\colon \CC K(x,y,z)_{t}\r\CC K(x,z)$, $t\geq 0$, in $\C V$, each of which is defined by applying the universal property of a push-out in $\C V$ to the commutative square
$$\xymatrix@C=150pt{
\bullet
\ar[r]^-{
\bar c_{\CC H}(b,y,z)\odot (f\odot \bar c_{\CC H}(b,y,a))^{\odot(t-1)}\odot f\odot\bar c_{\CC H}(x,y,a)
}
\ar[d]_{\psi(x,y,z)_{t}}&
\bullet
\ar[d]^{\bar\psi(x,z)_{t}
}\\
\CC K(x,y,z)_{t-1}\ar[r]_-{\Phi_{x,y,z,t-1}}&\CC K(x,z)}$$
The initial step is $\Phi_{x,y,z,0}=\xi_{x,y,z}$, in particular the equation $\xi_{x,y,z}=\Phi'_{x,y,z}\xi'_{x,y,z}$ follows. The commutativity of this square follows straightforwardly by induction, using the definition of composition in $\CC K$.

Similarly, $\Phi_{x,y,z}$ is defined by applying the universal property of a colimit to a sequence of compatible morphisms $\Phi'_{x,y,z,t}\colon \CC K(x,z)_{t}\r\CC K(x,y,z)$, $t\geq 0$, in $\C V$, each of which is defined by applying the universal property of a push-out in $\C V$ to the commutative square
$$\xymatrix@C=100pt{\bullet\ar[r]^-{\CC H(b,z)( f\otimes \CC H(b,a))^{\odot (t-1)}\odot f\otimes \CC H(x,a)}\ar[d]_{\psi(x,z)_{t}}\ar@{}[rd]|{\text{push}}&\bullet\ar[d]^{
\bar\psi(x,y,z)_{t}
}\\
\CC K(x,z)_{t-1}\ar[r]_-{\Phi'_{x,y,z,t-1}}&\CC K(x,y,z)}$$
The initial step $\Phi'_{x,y,z,0}$ is taken to be the composite $$\CC K(x,z)_{0}=\CC H(x,z)\st{\tilde\varphi}\To\CC K(x,y,z)_0\st{\xi'_{x,y,z}}\To\CC K(x,y,z).$$  Again,  commutativity  follows  by induction in a straightforward way. 

Equation $\xi'_{x,y,z}=\Phi_{x,y,z}\xi_{x,y,z}$ follows from the definition of $\Phi'_{x,y,z,0}$  and from the fact that the square
$$\xymatrix@C=70pt@R=10pt{\CC  K(y,z)_{s}\otimes\CC K(x,y)_{t}\ar[ddd]_-{c_{\CC K}(x,y,z)_{s,t}}\ar[r]_{
\begin{array}{c}\\[-15pt]
\text{\scriptsize canonical map}\\[-6pt]\text{\scriptsize to the colimit}
\end{array}}&\CC K(y,z)\otimes\CC K(x,y)\ar@{->>}[d]^{\text{quotient map to coequalizer}}\\
&\CC K(y,z)\otimes_{\CC K(y,y)}\CC K(x,y)\ar[d]^-{\tilde c}\\
&\CC K(x,y,z)_0\ar[d]^-{\xi'_{x,y,z}}\\
\CC  K(x,z)_{s+t}\ar[r]^-{\Phi'_{x,y,z,s+t}}&\CC K(x,y,z)}$$
commutes for all $s,t\geq 0$.

The equations $\Phi'_{x,y,z}\Phi_{x,y,z}=\id{\CC K(x,z)}$ and $\Phi_{x,y,z}\Phi'_{x,y,z}=\id{\CC K(x,y,z)}$ are consequences of the universal properties of colimits.
\end{proof}

The statement of Lemma \ref{composition} is coherent when the fact that the reduced composition morphism $\bar c_{\CC K}(x,y,z)$ is an isomorphism if $x=y$ or $y=z$. Indeed, in this case $\bar c_{\CC H}(x,y,z)$ is an isomorphism, hence the parallel arrow $\tilde c$ in \eqref{composdiag1} is an isomorphism. Moreover $\bar c_{\CC H}(x,y,a)$ and/or $\bar c_{\CC H}(b,y,z)$ is also an isomorphism, so $\xi_{x,y,z,t}$ in \eqref{composdiag2} is an isomorphism too. Therefore $\bar c_{\CC K}(x,y,z)$ is a transfinite composition of isomorphisms.

The following final lemma can be checked as the previous two ones.

\begin{lem}\label{rightmodule}
Consider the push-out \eqref{pussy}. The right $\CC K(x,x)$-module morphism $\CC H(x,y)\otimes_{\CC H(x,x)}\CC K(x,x)\r \CC K(x,y)$ induced by the right $\CC H(x,x)$-module morphism $\varphi(x,y)\colon \CC H(x,y)\r \CC K(x,y)$ is a transfinite composition in $\lmod{\CC K(x,x)^{\op}}$
$$\xymatrix@C=40pt{\CC H(x,y)\otimes_{\CC H(x,x)}\CC K(x,x)=\CC K(x,y)_{0}^{r}\r\cdots\r \CC K(x,y)_{t-1}^{r}\ar[r]^-{\varphi(x,y)_{t}^{r}} & \CC K(x,y)_{t}^{r}\r\cdots}$$
such that $\varphi(x,y)_{t}^{r}$ is a push-out  in $\lmod{\CC K(x,x)^{\op}}$ of the form
\begin{equation}\label{kaketa}
\xymatrix@C=170pt{\bullet\ar[r]^-{\bar c_{\CC H} (b,x,y)\odot( f\odot \bar c_{\CC H}(b,x,a))^{\odot (t-1)}\odot f\otimes \CC H(x,a)\otimes_{\CC H(x,x)}\CC K(x,x)}\ar[d]_{\psi(x,y)_{t}^{r}}\ar@{}[rd]|{\text{push}}&\bullet\ar[d]^{\bar\psi(x,y)_{t}^{r}}\\
\CC K(x,y)_{t-1}^{r}\ar[r]_-{\varphi(x,y)_{t}^{r}}&\CC K(x,y)_{t}^{r}}
\end{equation}
The target of the upper horizontal arrow is $\eqref{target2}\otimes_{\CC H(x,x)}\CC K(x,x)$. The morphism 
 $\psi(x,y)_{t}^{r}$ departs from the source of a push-out product, so it is defined by maps from the  right $\CC K(x,x)$-modules obtained from  $\eqref{target2}\otimes_{\CC H(x,x)}\CC K(x,x)$ by replacing the target of each push-out product factor with the source. If we replace a $V$ with a $U$, the map is defined as $\eqref{pegando}\otimes_{\CC H(x,x)}\CC K(x,x)$,  changing $\bar\psi(x,y)_{t-1}$
for $\bar\psi(x,y)_{t-1}^r$. If we replace the $s^{\text{th}}$ copy of $\CC H(b,a)$ with $\CC H(x,a)\otimes_{\CC H(x,x)} \CC H(b,x)$, $1\leq s\leq t-1$, then the map is
$$\xymatrix{
\CC H(b,y)\!\otimes\! V\!\otimes\! \cdots\!\otimes\! \CC H(x,a)\!\otimes_{\CC H(x,x)}\!\CC H(b,x)\!\otimes\!\cdots\!\otimes\! V\!\otimes\! \CC H(x,a)\!\otimes_{\CC H(x,x)}\!\CC K(x,x)
\ar[d]_{\id{}\otimes_{\CC H(x,x)}\bar\psi(x,x)_{t-s}\otimes_{\CC H(x,x)}\id{}}\\
\CC H(b,y)\otimes V\otimes \cdots\otimes \CC H(x,a)\otimes_{\CC H(x,x)}\CC K(x,x)_{t-s}\otimes_{\CC H(x,x)}\CC K(x,x)
\ar[d]_{\id{}\otimes_{\CC H(x,x)}\text{canonical map to the colimit }\otimes_{\CC H(x,x)}\id{}}\\
\CC H(b,y)\otimes V\otimes \cdots\otimes \CC H(x,a)\otimes_{\CC H(x,x)}\CC K(x,x)\otimes_{\CC H(x,x)}\CC K(x,x)
\ar[d]_{\id{}\otimes\text{composition}}\\
\CC H(b,y)\otimes V\otimes \cdots\otimes \CC H(x,a)\otimes_{\CC H(x,x)}\CC K(x,x)
\ar[d]_{\bar\psi(x,y)_{s}^{r}}\\
\CC K(x,y)_{s}^{r}
\ar[d]_{
\text{bonding}
}^{\text{or the identity if }s=t-1}\\
\CC K(x,y)_{t-1}^{r}
}$$
If we finally replace $\CC H(b,y)$ with $\CC H(x,y)\otimes_{\CC H(x,x)}\CC H(b,x)$ then the map is
similarly defined, putting $s=0$ and changing $\bar\psi(x,y)_{s}^{r}$ for the identity. 
%
\end{lem}

There is an obvious version of this lemma for left $\CC K(y,y)$-modules. It is actually a formal consequence of this one, taking opposite $\C V$-categories.

Let us analize some special cases of Lemma \ref{rightmodule}. If $a=x$ then the reduced composition $\bar c_{\CC H}(a,x,b)$ is an isomorphism, so $\varphi(x,y)_t^r$ is an isomorphism for $t>1$. Therefore we have a single push-out of right $\CC K(x,x)$-modules
$$
\xymatrix@C=60pt{\bullet 
\ar[d]_{\psi(x,y)_1^r}
\ar[r]^{\bar c_{\CC H} (b,x,y)\odot f\otimes \CC H(x,a)\otimes_{\CC H(x,x)}\CC K(x,x)}&\bullet
\ar[d]^{\bar\psi(x,y)_1^r}
\\
\CC H(x,y)\otimes_{\CC H(x,x)}\CC K(x,x)\ar[r]& \CC K(x,y)}$$
If $b=x$ or $x=y$ then $\bar c_{\CC H}(b,x,y)$ is an isomorphism, so all $\varphi(x,y)_t^r$ are isomorphisms, $t\geq 1$, and hence $\CC K(x,y)=\CC H(x,y)\otimes_{\CC H(x,x)}\CC K(x,x)$.

\section{Pseudo-cofibrations}

In this section we introduce a class of morphisms in a symmetric monoidal model category which share some relevant features with cofibrations. We also show that this class of morphisms is closed under certain constructions.


\begin{defn}\label{pcof}
A morphism $f$ in $\C V$ is a \emph{pseudo-cofibration} if $f\odot g$ is a (trivial) cofibration for any (trivial) cofibration $g$.
\end{defn}

\begin{rem}\label{honorary}
Cofibrations are pseudo-cofibrations by the push-out product axiom, but also the morphism $\varnothing\r\unit$, which need not be a cofibration, is a pseudo-cofibration since $(\varnothing\r \unit)\odot g=g$. The monoid axiom shows that, in general, $\varnothing \r V$ is a pseudo-cofibration if and only if $V$ is pseudo-cofibrant.
\end{rem}

\begin{lem}\label{retract}
Pseudo-cofibrations are closed under retracts.
\end{lem}

This follows from the fact that (trivial) cofibrations are closed under retracts.

\begin{lem}\label{pcofpush}
Given a push-out diagram in $\C V$,
$$\xymatrix{U\ar[d]_h\ar[r]^f\ar@{}[rd]|{\text{push}}&V\ar[d]^{h'}\\
X\ar[r]_{f'}&Y}$$
if $f$ is a pseudo-cofibration then so is $f'$.
\end{lem}

\begin{proof}
This follows from the fact that, by Lemma \ref{podot}, for any $g\colon U'\r V'$ in $\C V$,
$$\xy
(0,-1.8)*{U\otimes V'\!\!\!\!\bigcup\limits_{U\otimes U'}\!\!\!\!V\otimes U'},
(0,-21.8)*{X\otimes V'\!\!\!\!\bigcup\limits_{X\otimes U'}\!\!\!\!Y\otimes U'},
(40,0)*{V\otimes V'},
(40,-20)*{Y\otimes V'},
(20,-10)*{\text{\scriptsize push}}
\ar(14,0);(33,0)^-{f\odot g}
\ar(40,-3);(40,-17)^-{h'\otimes\id{V'}}
\ar(14,-20);(33,-20)_-{f'\odot g}
\ar(0,-6);(00,-17)_-{h\otimes\id{V'}\!\!\!\!\bigcup\limits_{h\otimes\id{U'}}\!\!\!\! h'\otimes\id{U'}}
\endxy$$
is also a push-out, so the upper (and hence the lower) horizontal arrow is a (trivial) cofibration if $g$ is so.
\end{proof}

\begin{lem}\label{continuo}
Let $\alpha$ be an ordinal, $F,G\colon\alpha\r\C V$ two continuous (i.e.~colimit preserving) diagrams indexed by $\alpha$, and $\tau\colon F\rr G$ a natural transformation such that $\tau(0)\colon F(0)\r G(0)$ and 
$$(\tau(\lambda+1),G(\lambda<\lambda+1))\colon F(\lambda+1)\bigcup_{F(\lambda)}G(\lambda)\To G(\lambda+1)$$ are pseudo-cofibrations, $\lambda+1<\alpha$. Then $\colim_\alpha \tau$ is a pseudo-cofibration.
\end{lem}

\begin{proof}
Let $g\colon U'\r V'$ be a (trivial) cofibration in $\C V$. We are going to consider the natural transformation $\tau\odot g$. The source and the target of $\tau\odot g$ are continuous since $F$ and $G$ are. The natural transformation $\tau\odot g$ looks as follows at a successor ordinal
 $$\xymatrix@C=60pt{   F(\lambda)\otimes V'\bigcup_{F(\lambda)\otimes U'}G(\lambda)\otimes U'\ar[d]^-{F(\lambda<\lambda+1)\otimes\id{}\bigcup G(\lambda<\lambda+1)\otimes\id{}}
\ar[r]^-{\tau(\lambda)\odot g}
&G(\lambda)\otimes V'  \ar[d]^-{G(\lambda<\lambda+1)\otimes \id{}}\\
   F(\lambda+1)\otimes V'\bigcup_{F(\lambda+1)\otimes U'}G(\lambda+1)\otimes U'\ar[r]_-{\tau(\lambda+1)\odot g}&G(\lambda+1)\otimes V'  
}$$
The morphism induced by this commutative square from the push-out of its left upper corner to the right bottom object is 
$(\tau(\lambda+1),G(\lambda<\lambda+1))\odot g$, 
which is a (trivial) cofibration since $(\tau(\lambda+1),G(\lambda<\lambda+1))$ is a pseudo-cofibration. Moreover $\tau(0)\odot g$ is a (trivial) cofibration since $\tau(0)$ is a pseudo-cofibration. This shows that $\tau\odot g$ is a (trivial) cofibration in the Reedy model category $\C V^\alpha$ of $\alpha$-sequences of objects in $\C V$. Hence its colimit $\colim_\alpha  (\tau \odot g)=(\colim_\alpha  \tau )\odot g$ is again a (trivial) cofibration, since the colimit is a left Quillen functor, see \cite[Corollary 5.1.6]{hmc}.
\end{proof}

\begin{lem}\label{pcofcompos}
The composition of pseudo-cofibrations is a pseudo-cofibration.
\end{lem}

\begin{proof}
Let $U\st{f}\to V\st{f'}\to W$
be pseudo-cofibrations and $g\colon X\r Y$ a (trivial) cofibration. In particular, $f\odot g$ and $f'\odot g$ are (trivial) cofibrations. The push-out product $(f'f)\odot g$ decomposes as
\begin{equation}\label{unacom}
\xy
(-40,-1.5)*{U\otimes Y\!\!\!\!\bigcup\limits_{U\otimes X}\!\!\!\!W\otimes X},
(10,-1.5)*{V\otimes Y\!\!\!\!\bigcup\limits_{V\otimes X}\!\!\!\!W\otimes X},
(40,0)*{W\otimes Y},
\ar(24,0);(33,0)^-{f'\odot g}
\ar(-26,0);(-4,0)
\endxy
\end{equation}
where the first arrow is obtained by taking the push-out of the rows of the following commutative diagram
$$\xymatrix@C=40pt{U\otimes Y\ar[d]_{f\otimes\id{}}
&U\otimes X\ar[l]_-{\id{}\otimes g}\ar[r]^-{(f'f)\otimes\id{}}\ar[d]_{f\otimes\id{}}&W\otimes X\ar@{=}[d]\\
V\otimes Y&V\otimes X\ar[l]^-{\id{}\otimes g}\ar[r]_-{f'\otimes\id{}}&W\otimes X}$$
The square on the left induces the morphism $f\odot g$  from the push-out of its right upper corner to its bottom left object. This morphism is a (trivial) cofibration, therefore the first arrow in \eqref{unacom} is a (trivial) cofibration by (the version for trivial cofibrations of) \cite[Lemma 7.2.15]{hirschhorn}. The morphism $f'\odot g$ is also a (trivial) cofibration, hence we are done.
\end{proof}


\section{Endomorphism monoids of cofibrant enriched categories}

In this section we prove that cofibrant objects in $\cat{\C V}{S}$ have cofibrant endomorphism monoids, and moreover cofibrations between them induce cofibrations on endomorphism monoids. 
In order to achieve this goal,  we
generalize Berger--Moerdijk's interval cofibrancy theorem \cite[Theorem
1.15]{htec} using the categorical constructions in Section \ref{pushpush}.

In the introduction we recalled the notion of \emph{pseudo-cofibrant} object in $\C V$,  introduced in \cite{htnso2}. Pseudo-cofibrant objects are closed under retracts and tensor products. 
Cofibrant objects are pseudo-cofibrant but the tensor unit $\unit$, which in general is not cofibrant, is pseudo-cofibrant for obvious reasons. Moreover, if $U$ is pseudo-cofibrant and $U\into V$ is a cofibration then $V$ is also pseudo-cofibrant \cite[Lemma A.4]{htnso2}. 

By the monoid axiom, pseudo-cofibrant objects are exactly those $V$ such that $V\otimes-\colon\C V\r \C V$ is a left Quillen functor. Hence, $V$ is pseudo-cofibrant if and only if $V\otimes f$ is a cofibration for all $f\in I\cup J$.

\begin{defn}\label{ppax}
The \emph{push-out product axiom} says that the push-out product $f\odot g$ of two cofibrations $f$ and $g$ in $\C V$ is again a cofibration, and if moreover $f$ or $g$ is a trivial cofibration then so is $f\odot g$.
\end{defn}

Recall from \cite{ammmc} that, given a monoid $A$ in $\C V$, the category $\lmod{A}$ of left $A$-modules inherits a model structure with fibrations and weak equivalences defined as in $\C V$. It is combinatorial with sets of generating (trivial) cofibrations $A\otimes I$ and $A\otimes J$, in particular the forgetful functor $\lmod{A}\r\C V$ is a right Quillen functor with left adjoint $A\otimes -\colon\C V\r\lmod{A}$. We know that the forgetful functor is also a left adjoint. It is even a left Quillen functor in certain circumstances.

\begin{lem}\label{contrario}
If $A$ is a monoid which is pseudo-cofibrant as an object in $\C V$ then the forgetful functor $\lmod{A}\r\C V$ is a left Quillen functor.
\end{lem}

 It is enough to check that generating (trivial) cofibrations are preserved, and this is obvious since $A$ is pseudo-cofibrant in $\C V$. 
%
%

The following result looks stronger than the push-out product axiom, but it is a consequence of it.

\begin{lem}\label{ppgordo}
Let $A$ and $B$ be monoids in $\C V$. If $f$ is a cofibration in $\lmod{A}$ and $g$ is a cofibration in $\lmod{B}$ then the push-out product $f\odot g$ in $\C V$, which is naturally a left $A\otimes B$-bimodule morphism, is a cofibration in $\lmod{A\otimes B}$. Moreover, if either $f$ or $g$ is also a weak equivalence, then so is $f\odot g$. 
\end{lem}

This lemma can be proved as \cite[Lemma 1.6]{brew}. It can be applied when $A$ (or $B$) is the initial monoid $\unit$, so $f$ (or $g$) is a plain morphism in $\C V$ and $f\odot g$ is a left $B$-(or $A$-)module morphism. 

We now extend the notion of pseudo-cofibrant object to module categories.

\begin{defn}
Given a monoid $A$ in $\C V$, a left $A$-module $M$ is \emph{pseudo-cofibrant} if $M\otimes-\colon\C V\r\lmod{A}$ preserves cofibrations. 
\end{defn}

This notion of pseudo-cofibrant module is not a particular case of the previous one since module categories over non-commutative monoids need not be monoidal.

Again, by the monoid axiom, a left $A$-module $M$ is pseudo-cofibrant if and only if the functor $M\otimes-\colon\C V\r\lmod{A}$ is a left Quillen functor. Hence it is enough to check  that morphisms in $M\otimes (I\cup J)$ are cofibrations in $\lmod{A}$. In particular $A$ is always pseudo-cofibrant as a left $A$-module, although it need not be cofibrant if $\unit$ is not cofibrant in $\C V$. 

One can check as in \cite[Lemma A.4]{htnso2}, using Lemma \ref{ppgordo} instead of the push-out product axiom, that a cofibration in $\lmod{A}$ with pseudo-cofibrant source has also a pseudo-cofibrant target.

Pseudo-cofibrant modules are preserved by extensions of scalars.

\begin{lem}\label{ee}
If $M$ is a pseudo-cofibrant left $A$-module and $A\r B$ is a monoid morphism then $B\otimes _AM$ is a pseudo-cofibrant left $B$-module.
\end{lem}

\begin{proof}
Notice that $B\otimes_AM\otimes -\colon\C V\r\lmod{B}$ is a composition of two left Quillen functors,
$$\xymatrix{\C V\ar[r]^-{M\otimes -}&\lmod{A}\ar[r]^-{B \otimes_A-}&\lmod{B}.}$$
\end{proof}

\begin{defn}\label{slpc}
A $\C V$-category $\CC H$ is \emph{strongly locally pseudo-cofibrant} if: 
\begin{enumerate}
 \item any morphism object $\CC H(x,y)$ is pseudo-cofibrant as a left $\CC H(y,y)$-module, as a right $\CC H(x,x)$-module, and as a plain object in $\C V$, $x,y\in\ob\CC H$,
 \item and reduced compositions morphisms in $\CC H$ are pseudo-cofibrations in $\C V$.
\end{enumerate}
\end{defn}

The initial object $\unit_{S}$ in $\cat{\C V}{S}$ clearly satisfies this property. Indeed, all endomorphism monoids are $\unit$ and the rest of morphism objects are $\varnothing$. Reduced compositions are either one of the identities $\varnothing\r \varnothing$ and $\unit\r\unit$ or the map $\varnothing\r \unit$.

Following the standard terminology, a $\C V$-category with pseudo-cofibrant morphism objects in $\C V$ should be plainly called \emph{locally pseudo-cofibrant}. The previous definition is a strengthening of this notion. 

\begin{prop}\label{paralpc}
A cofibration $\varphi\colon\CC H\into\CC K$ in $\cat{\C V}{S}$ with locally pseudo-cofibrant source $\CC H$ is a local  cofibration, in particular $\CC K$ is also locally pseudo-cofibrant.
\end{prop}

\begin{proof}
Cofibrations in $\C V$ are closed under retracts and transfinite compositions, hence it is enough to assume that $\varphi$ is a push-out of a generating cofibration in $\cat{\C V}{S}$, i.e.~like \eqref{pussy} with $f\in I$. The upper horizontal arrow in \eqref{cosaque} is isomorphic to
\begin{equation}\label{comer}
(\CC H(b,y)\otimes\CC H(b,a)^{\otimes(t-1)}\otimes \CC H(x,a))\otimes (f^{\odot t}). 
\end{equation}
The object on the left is pseudo-cofibrant since $\CC H$ is locally pseudo-cofibrant and the map on the right is a cofibration by the push-out product axiom, therefore the previous map is a cofibration in $\C V$. This shows that $\varphi$ is locally a transfinite composition of cofibrations, hence we are done.
\end{proof}

There is some redundancy in Definition \ref{slpc}, e.g.~endomorphism objects, as all monoids, are pseudo-cofibrant regarded as left/right modules over themselves. The following lemma exhibits another redundancy.

\begin{lem}
If $A$ is a monoid which is pseudo-cofibrant  as an object in $\C V$ and $M$ is a pseudo-cofibrant left $A$-module then $M$ is also pseudo-cofibrant in $\C V$.
\end{lem}

\begin{proof}
Observe that $M\otimes -\colon\C V\r\C V$ is the composition of two left Quillen functors
$$\xymatrix{\C V\ar[r]^-{M\otimes -}&\lmod{A}\ar[r]^-{\text{forget}}&\C V,}$$
see Lemma \ref{contrario}.
\end{proof}

The following additional property of pseudo-cofibrant modules is needed to show that they are closed under tensor products.

\begin{lem}
If $M$ is a pseudo-cofibrant left $A$-module then $M\otimes-\colon\lmod{B}\r\lmod{A\otimes B}$ is a left Quillen functor for any monoid $B$.
\end{lem}

\begin{proof}
If $f\in I$ (resp.~$J$) then $M\otimes f$ is a (trivial) cofibration in $\lmod{A}$ since $M$ is a pseudo-cofibrant left $A$-module. 
The extension of scalars $\iota_*\colon\lmod{A}\r \lmod{A\otimes B}$ along the monoid morphism 
$\iota\colon A\r A\otimes B$ induced by the unit of $B$ is a left Quillen functor (as all extensions of scalars),  therefore $M\otimes (B\otimes f)=\iota_*(M\otimes f)$ is a (trivial) cofibration in $\lmod{A\otimes B}$. This proves that the functor in the statement takes generating (trivial) cofibrations to (trivial) cofibrations, as desired.
\end{proof}

\begin{cor}\label{tensorpc}
If $M$ is a pseudo-cofibrant left $A$-module and $N$ is a pseudo-cofibrant left $B$-module then $M\otimes N$ is a pseudo-cofibrant left $A\otimes B$-bimodule.
\end{cor}

\begin{proof}
Simply notice that $M\otimes N\otimes -\colon\C V\r \lmod{A\otimes B}$ is a composition of two left Quillen functors
$$\xymatrix{\C V\ar[r]^-{N\otimes -}&\lmod{B}\ar[r]^-{M\otimes -}&\lmod{A\otimes B}.}$$
\end{proof}

\begin{defn}
A morphism $\varphi\colon\CC H\r\CC K$ in $\cat{\C V}{S}$ is a \emph{strong local cofibration} if, for any $x,y\in S$, the morphism $\varphi(x,y)\colon\CC H(x,y)\r\CC K(x,y)$  is a cofibration in $\C V$ and 
the induced morphisms $$\CC H(x,y)\otimes_{\CC H(x,x)}\otimes \CC K(x,x)\To\CC K(x,y),\qquad 
\CC K(y,y)\otimes_{\CC H(y,y)}\otimes \CC H(x,y)\To\CC K(x,y),$$ are cofibrations in $\lmod{\CC K(x,x)^{\op}}$ and $\lmod{\CC K(y,y)}$, respectively.
\end{defn}

Strong local cofibrations are also local cofibrations in the usual sense. 
Notice that a retract or (transfinite) composition of strong local cofibrations is again a strong local cofibration.

\begin{thm}\label{ict2}
Let $\varphi\colon\CC H\into \CC K$ be a cofibration  in $\cat{\C V}{S}$ with strongly locally pseudo-cofibrant source $\CC H$. Then $\varphi(x,x)\colon\CC H(x,x)\r \CC K(x,x)$ is a cofibration in $\mon{\C V}$ for all $x\in S$, $\varphi$ is a strong local cofibration, and $\CC K$ is strongly locally pseudo-cofibrant.
\end{thm}

\begin{proof}
Cofibrations in $\cat{\C V}{S}$ are retracts of relative $T_S(I_S)$-cell complexes. Since (pseudo-)cofibrations and pseudo-cofibrant objects are closed under retracts, it is enough to assume that $\varphi$ is a relative $T_S(I_S)$-cell complex. The proof is by  induction on the length of the transfinite composition. 

We start with the induction setp at a successor ordinal, i.e.~we suppose that $\varphi$ is a push-out of a generating cofibration, like \eqref{pussy} with $f\in I$. We first check that $\varphi(x,x)$ is a monoid cofibration. 

The category of monoids in ${\lmod{\CC H(x,x)^{\env}}}$ is regarded as a model category via the isomorphism $\mon{\lmod{\CC H(x,x)^{\env}}}\cong \CC H(x,x)\downarrow\mon{\C V}$ in the paragraph preceding Lemma \ref{monoid}. Fibrations and weak equivalences are defined as in $\C V$, therefore $T\colon{\lmod{\CC H(x,x)^{\env}}}\rightleftarrows \mon{\lmod{\CC H(x,x)^{\env}}}\colon\text{forg.}$ is a Quillen pair.


Since $\CC H$ is strongly locally pseudo-cofibrant,  $\bar c_{\CC H}(b,x,a)$ is a pseudo-cofibration in $\C V$, $\CC H(b,x)\otimes\CC H(x,a)$ is a pseudo-cofibrant $\CC H(x,x)$-bimodule by Corollary \ref{tensorpc}, and therefore the upper horizontal arrow in \eqref{monoiddiag} is a free monoid morphism on a cofibration in $\lmod{\CC H(x,x)^{\env}}$, hence a cofibration in $\mon{\lmod{\CC H(x,x)^{\env}}}$. This and  Lemma \ref{monoid} show that $\varphi(x,x)$ is a transfinite composition of cofibrations in $\mon{\lmod{\CC H(x,x)^{\env}}}$. In particular it is  a cofibration in $\CC H(x,x)\downarrow\mon{\C V}$, and hence also in $\mon{\C V}$.

Let us see that $\varphi$ is a strong local cofibration. We have already checked in Proposition \ref{paralpc} that it is a local cofibration. 
Since $\CC H$ is strongly locally pseudo-cofibrant, $\bar c_{\CC H}(b,x,a)$ and $\bar c_{\CC H}(b,x,y)$ are pseudo-cofibrations in $\C V$, $\CC H(x,a)\otimes_{\CC H(x,x)}\CC K(x,x)$ is a pseudo-cofibrant right $\CC K(x,x)$-module by Lemma \ref{ee}, and therefore the upper horizontal arrow in \eqref{kaketa} is a cofibration in $\lmod{\CC K(x,x)^{\op}}$. By Lemma \ref{rightmodule}, $\CC H(x,y)\otimes_{\CC H(x,x)}\otimes \CC K(x,x)\r\CC K(x,y)$ is a transfinite composition of cofibrations in $\lmod{\CC K(x,x)^{\op}}$, so it is also a cofibration. One can similarly check that $\CC K(y,y)\otimes_{\CC H(y,y)}\otimes \CC H(x,y)\r\CC K(x,y)$ is a cofibration of left $\CC K(y,y)$-modules invoking the left module version of Lemma \ref{rightmodule}, which is a formal consequence of it.

Since $\varphi$  is a strong local cofibration, the first property of strongly locally pseudo-cofibrant $\C V$-categories for $\CC K$ is derived from the same property for $\CC H$ and Lemma \ref{ee}. Let us check the second one. 

The upper horizontal arrow in \eqref{composdiag1} is a pseudo-cofibration in $\C V$ since $\CC H$ is strongly locally pseudo-cofibrant, hence the parallel arrow $\tilde c$ is also a pseudo-cofibration by Lemma \ref{pcofpush}. Moreover, for the same reason the reduced composition morphisms in the upper horizontal arrow in \eqref{composdiag2} are pseudo-cofibrations, so this arrow is an honest cofibration. Therefore $\xi_{x,y,z}$ in \eqref{composdiag1} is a cofibration, since it is a transfinite composition of cofibrations by Lemma \ref{composition}. In particular, $\bar c_{\CC K}(x,y,z)=\xi_{x,y,z}\tilde c$ is a pseudo-cofibration by Lemma \ref{pcofcompos}. 

We now prove the induction step at a limit ordinal. Assume that $\varphi$ is a transfinite composition $\varphi\colon\CC H=\CC K_0\r\colim_{\lambda<\alpha}\CC K_\lambda=\CC K$, $\alpha$ a limit  ordinal, where each  push-out of a generating cofibration  $\CC K_\lambda\r\CC K_{\lambda+1}$  is a strong local cofibration with strongly pseudo-cofibrant source and target which induces cofibrations $\CC K_\lambda(x,x)\r\CC K_{\lambda+1}(x,x)$ on endomorphism monoids in $\C V$. Then $\varphi$ is a strong local cofibration and $\varphi(x,x)$ is a cofibration in $\mon{\C V}$ 
since these kinds of maps are closed under transfinite composition. As we have already seen, this implies that $\CC K$ satisfies the first property of strongly locally pseudo-cofibrant $\C V$-categories. Let us check that Lemma \ref{continuo} applies to prove that $\bar c_{\CC K}(x,y,z)=\colim_{\lambda <\alpha}\bar c_{\CC K_\lambda}(x,y,z)$ is a pseudo-cofibration. The initial term is the pseudo-cofibration $\bar c_{\CC H}(x,y,z)$. We also have to check that for each $\lambda<\alpha$ the commutative square 
 $$\xymatrix@C=60pt{   \CC K_{\lambda}(y,z)\otimes_{\CC K_{\lambda}(y,y)}\CC K_{\lambda}(x,y)\ar[d]\ar[r]^-{\bar c_{\CC K_{\lambda}}(x,y,z)}&  \CC K_{\lambda}(x,z)\ar[d]\\
   \CC K_{\lambda+1}(y,z)\otimes_{\CC K_{\lambda+1}(y,y)}\CC K_{\lambda+1}(x,y)\ar[r]_-{\bar c_{\CC K_{\lambda+1}}(x,y,z)}&\CC K_{\lambda+1}(x,z)  }$$
induces a pseudo-cofibration from the push-out of the left upper corner to the right bottom object. This induced map is actually like $\xi_{x,y,z}$ in Lemma \ref{composition}, so it is even a cofibration by the argument in the previous paragraph. 
\end{proof}

\begin{cor}\label{ict1}
Cofibrant objects in $\cat{\C V}{S}$
are strongly locally pseudo-cofibrant and have cofibrant endomorphism monoids.
\end{cor}

%
%

%
%

\section{The monoid axiom}

In this short section we recall the monoid axiom introduced in \cite[Definition 3.3]{ammmc} and prove some related results.

\begin{defn}\label{monax}
Consider the following class of morphisms in $\C V$,
$$K=\{X\otimes f\,;\,f\text{ is a trivial cofibration and }X\text{ is an object in }\C V\}.$$
The \emph{monoid axiom} for $\C V$ says that relative $K$-cell complexes are weak equivalences.
\end{defn}

The monoid axiom implies that retracts of relative $K$-cell complexes are weak equivalences. However it does not directly imply that transfinite compositions of retracts of relative $K$-cell complexes are weak equivalences, which is a technicality that we will later need. If $K$ were a set, we would deduce that $\cellr{K}=\cof{K}$, so $\cellr{K}$ would be closed under  transfinite compositions, which is exactly what we want. The class $K$ need not be a set, but the following lemma shows that $\cellr{K}$  is still closed under transfinite compositions.

\begin{lem}
Let $C$ be any class of morphisms in a locally presentable category. The class $\cellr{C}$ of retracts of relative $C$-cell complexes is closed under transfinite compositions.
\end{lem}

\begin{proof}
Suppose $f\colon X\r Y$ is a transfinite composition of morphisms in $\cellr{C}$. Each step of this transfinite composition is a retract of a relative $C$-cell complex. Each of these relative $C$-cell complexes is a transfinite composition of push-outs of morphisms in a certain subset of $C$. If $C'$ is the union of all these subsets of $C$, then we notice that $f$ is in fact a transfinite composition of retracts of relative $C'$-cell complexes. Since $C'$ is a set, $\cellr{C'}=\cof{C'}$, which is closed under transfinite compositions, hence $f\in\cellr{C'}\subset\cellr{C}$.
\end{proof}

\begin{defn}\label{localin}
A $\C V$-functor $\varphi\colon\CC H\r\CC K$ is \emph{locally in $\cellr{K}$} if $\varphi(x,y)\colon\CC H(x,y)\r\CC K(\varphi(x),\varphi(y))$ is in $\cellr{K}$ for all $x,y\in\ob\CC H$.
\end{defn}

By the monoid axiom, $\C V$-functors which are locally in $\cellr{K}$ are homotopically fully faithful. The following result is a consequence of the previous lemma.

\begin{cor}\label{closcor}
A transfinite composition of $\C V$-functors which are locally in $\cellr{K}$ is again  locally in $\cellr{K}$, in particular it is homotopically fully faithful.
\end{cor}

We now identify an important class of $\C V$-functors which are locally in $\cellr{K}$.

\begin{prop}\label{yanose}
Any trivial cofibration $\varphi\colon\CC H\r\CC K$ in $\cat{\C V}{S}$ is locally in $\cellr{K}$.
\end{prop}

\begin{proof}
We can suppose without loss of generality that $\varphi$ is a relative $T_S(J_S)$-cell complex. Since $\C V$-functors which are locally in $\cellr{K}$ are closed under transfinite compositions and retracts, we can even assume that $\varphi$ fits into a push-out of the form \eqref{pussy} with $f\in J$. The upper horizontal arrow in \eqref{cosaque} is in $K$, since it is isomorphic to \eqref{comer}
and $f^{\odot t}$ is a trivial cofibration by the push-out product axiom. Hence the description of the push-out \eqref{pussy} given in Section \ref{pushpush} shows that $\varphi(x,y)$ is a relative $K$-cell complex for all $x,y\in S$.
\end{proof}


\section{Closure properties of homotopically essentially surjective enriched functors}\label{closure}

The class of DK-equivalences is obviously closed under compositions and retracts, since both equivalences of categories and weak equivalences in $\C V$ are closed under these operations. Homotopically essentially surjective $\C V$-functors are however closed under more general operations. 

\begin{prop}\label{closure1}
Consider a push-out square in $\cat{\C V}{}$ as follows
$$\xymatrix{
\CC H\ar[r]^-{\phi}\ar[d]_{\psi}\ar@{}[rd]|{\text{push}}&\CC K\ar[d]^{\psi'}\\
\CC L\ar[r]_-{\phi'}&\CC M
}$$
If $\phi$ is homotopically essentially surjective then so is $\phi'$.
\end{prop}

\begin{proof}
Let $x\in\ob \CC M$. The underlying diagram of object sets is a push-out, so if $x$ is not in the image of $\phi'$ then it is $x=\psi'(x')$ for some $x'\in\ob\CC K$. Since $\phi$ is homotopically essentially surjective we can take $y\in\ob\CC H$ with $\phi(y)\cong x'$ in $\pi_0\CC K$. If we take $\pi_0$ on the push-out square we still get a commutative square of ordinary categories hence $\phi'\psi(y)=\psi'\phi(y)\cong\psi'(x')=x$ in $\pi_0\CC M$.
\end{proof}

%

\begin{prop}\label{closure2}
A transfinite composition of homotopically essentially surjective $\C V$-functors is homotopically essentially surjective.
\end{prop}

\begin{proof}
This follows by induction. Since homotopically essentially surjective $\C V$-functors are obviously closed under compositions, it is enough to check the induction step at limit ordinals. 

Let $\varphi\colon \CC H=\CC K_0\r\colim_{\lambda<\alpha}\CC K_\lambda$ be a transfinite composition, $\alpha$ a limit ordinal, such that 
the smaller transfinite compositions $\CC H=\CC K_0\r\CC K_{\lambda}$, $\lambda <\alpha$, are homotopically essentially surjective. Since $\ob\CC K=\colim_{\lambda<\alpha}\ob\CC K_\lambda$, any object in $\CC K$ comes from some $\CC K_\lambda$, $\lambda<\alpha$, hence by hypothesis it is isomorphic to an object coming from $\CC H$.
\end{proof}

\section{The Dwyer--Kan model structure}\label{laprueba}

In this section we prove our main result, Theorem \ref{main}. The following proposition is an important previous step. Ideally, we would have proved it in Section \ref{gtc}, but we actually need the strongest results in the intermediate sections.

\begin{prop}\label{tress}
Any relative $J'$-cell complex $\varphi\colon\CC H\r \CC K$ is locally in $\cellr{K}$, in the sense of Definition \ref{localin}. In particular $\varphi$ is homotopically fully faithful.
\end{prop}

\begin{proof}
By Corollary \ref{closcor}, we can suppose that $\varphi$ is a push-out of a morphism in $J'$. Proposition \ref{yanose} covers the case when $\varphi$ is a push-out of a morphism in $T_{\{0,1\}}(J_{01})$, see Proposition \ref{01}. Therefore, we can suppose that $\varphi$ fits into a push-out in $\cat{\C V}{}$ as follows,
$$\xymatrix{
\widetilde{i^{*}\CC I}\ar[r]^-{\theta_{\CC I}}\ar[d]\ar@{}[rd]|{\text{push}}&\tilde{\CC I}\ar[d]\\
\CC H\ar[r]_-{\varphi}&\CC K
}$$
This push-out can be decomposed in two steps
$$\xymatrix@C=40pt{
\widetilde{i^{*}\CC I}\ar[r]_-\sim^-{\theta_{\CC I}(0,0)}\ar[d]\ar@{}[rd]|{\text{push}}&\tilde{\CC I}(0,0)\ar[r]^-{\text{full incl.}}\ar[d]\ar@{}[rd]|{\text{push}}&\tilde{\CC I}\ar[d]\\
\CC H\ar[r]_-{\varphi_{1}}&\tilde{\CC H}\ar[r]_-{\varphi_{2}}&\CC K
}$$
The $\C V$-functor $\varphi_{2}$ is fully faithful by \cite[Proposition 3.1]{cmcpfo}, compare also \cite[Proposition 5.2]{hin}. Recall that $\theta_{\CC I}(0,0)$ is a weak equivalence by construction. Moreover, this monoid morphism is part of a cofibration with cofibrant source in the fiber $\cat{\C V}{\{0,1\}}$, hence it is a cofibration of monoids by Theorem \ref{ict2} and Corollary \ref{ict1}. 
Therefore $\varphi_1$ is a trivial cofibration in $\cat{\C V}{\ob\CC H}$, in particular it is 
locally in $\cellr{K}$ by Proposition \ref{yanose}. Hence $\varphi=\varphi_2\varphi_1$ is 
locally in $\cellr{K}$.
\end{proof}

\begin{prop}\label{doss}
 Any relative $J'$-cell complex is an $I'$-cofibration and a DK-equivalence.
\end{prop}

\begin{proof}
By Proposition \ref{011}, morphisms in $T_{\{0,1\}}(J_{01})$ are $T_{\{0,1\}}(I_{01})$-cofibrations, hence $I'$-cofibra\-tions. Moreover, morphisms in $\{\theta_{\CC I}\,;\,\CC I\in G\}$ are also $I'$-cofibrations by Proposition \ref{thetai}. Therefore relative $J'$-cell complexes are $I'$-cofibrations. They are homotopically essentially surjective and homotopically fully faithful by Propositions \ref{closure1}, \ref{closure2}, and \ref{tress}.
\end{proof}

We can finally prove our main result.

\begin{proof}[Proof of Theorem \ref{main}]
We are going to use the recognition theorem for cofibrantly generated model categories in \cite[Theorem 2.1.19]{hmc} in order to show that $\cat{\C V}{}$  is a cofibrantly generated model category with set of generating cofibrations $I'$ in \eqref{iprima} and set of generating trivial cofibrations $J'$ in \eqref{jprima}, for $G$ a generating set of $\C V$-intervals, see Definition \ref{vinterval} and Proposition \ref{lohay}.

The category $\cat{\C V}{}$ is locally presentable by \cite{vcatlplb}. Hence it is (co)complete and the second and third conditions of the recognition theorem hold. The first one says that DK-equivalences must be closed under retracts and satisfy the 2-out-of-3 property. This is obvious since weak equivalences in $\C V$ and equivalences of categories satisfy these properties. The last three conditions have been verified in Propositions \ref{doss} and \ref{unno}.
%
%
%
\end{proof}

\begin{rem}\label{comoson}
Cofibrations in $\cat{\C V}{}$ are $I'$-cofibrations, which can be understood via Corollary \ref{008}. A $\C V$-category $\CC H$ is cofibrant if and only if it is cofibrant in $\cat{\C V}{\ob\CC H}$. In this case, endomorphism objects  $\CC H(x,x)$ are cofibrant as monoids and pseudo-cofibrant in $\C V$ (since the unit $\unit\r \CC H(x,x)$ is a cofibration) and the rest of morphism objects $\CC H(x,y)$, $x\neq y$, are cofibrant in $\C V$, see Theorem \ref{ict2} and Corollary \ref{ict1}.

All fibrations in $\cat{\C V}{}$ are local fibrations by Proposition \ref{01} (2), but this is not a characterization. Fibrations are more difficult to understand than cofibrations. They are just $J'$-injective maps. This property can be phrased as a kind of enriched isomorphism lifting property using $\C V$-intervals. This constrasts with $\cat{\C V}{S}$, where fibrations are precisely local fibrations, and  prevents us from finding easy conditions under which $\cat{\C V}{S}$ is right proper.

Something similar happens with fibrant objects $\CC H$. They have fibrant morphism objects $\CC H(x,y)$ in $\C V$. This is enough to be fibrant in $\cat{\C V}{\ob\CC H}$, but in order to be fibrant in $\cat{\C V}{}$, $\CC H$ must also satisfy an injectivity property with respect to $\C V$-intervals.
\end{rem}

\section{Change of base category}


Let $\C W$ be a category with the same structure and satisfying the same properties as $\C V$ and $F\colon\C V\rightleftarrows\C W\colon G$ a Quillen pair with lax symmetric monoidal right adjoint $G$. The derived adjoint pair
\begin{equation}\label{apderived}
\xymatrix{ \ho{\C V} \ar@<.5ex>[r]^-{\mathcal{L}F}& \ho{\C W} \ar@<.5ex>[l]^-{\mathcal{R}G}}
\end{equation}
satisfies the same formal properties as $F\colon\C V\rightleftarrows\C W\colon G$, concerning the symmetric monoidal structures. In particular, $\mathcal RG$ is lax symmetric monoidal. This has some nice consequences that we now describe. One the one hand, there is an induced functor between categories of small enriched categories
$$\mathcal RG\colon \cat{\ho\C W}{}\To \cat{\ho\C V}{}$$
like \eqref{kaker}.

\begin{prop}\label{conputa0}
Consider
the following diagram of categories and functors,
$$\xymatrix{\cat{\C V}{}\ar@{<-}[r]^-{G}\ar[d]_{p_{\C V}}&\cat{\C W}{}\ar[d]^{p_{\C W}}\\
\cat{\ho\C V}{}\ar@{<-}[r]_{\mathcal RG}&\cat{\ho\C W}{}}$$
where the vertical arrows are given by the obvious functors from a model category to its homotopy category.
If we restrict to fibrant objects in $\cat{\C W}{}$, then the square commutes up to natural isomorphism.
\end{prop}

This follows easily from the fact that fibrant $\C W$-categories have fibrant morphism objects, see Remark \ref{comoson}.

On the other hand, there is a monoidal natural transformation
$$\xymatrix{\ho\C W\ar[rr]^{\mathcal RG}\ar[rd]_-{\pi_0}&&\ho\C V\ar[ld]^-{\pi_0}\\
&\operatorname{Set}&\ar@{=>}(15,-5);(20,-5)^-\nu}$$
defined by
$$\pi_0W=\ho\C W(\unit,W)\st{\mathcal RG}\To\ho\C V(\mathcal RG(\unit),\mathcal RG(W))\r\ho\C V(\unit,\mathcal RG(W))=\pi_0\mathcal RG(W)$$
Here the second arrow is induced by the morphism $\unit\r\mathcal RG(\unit)$ which is part of the lax symmetric monoidal structure of $\mathcal RG$. This natural transformation $\nu$ is obviosuly an isomorphism when $F\colon\C V\rightleftarrows\C W\colon G$ is a Quillen equivalence.

\begin{cor}\label{preserba0}
There is a natural functor $\pi_0\CC K\r\pi_0\mathcal RG(\CC K)$ for any $\C W$-category $\CC K$. This functor is the identity on objects. Moreover, it is an isomorphism if $F\colon\C V\rightleftarrows\C W\colon G$ is a Quillen equivalence. 
\end{cor}

Suppose from now on that $F\colon\C V\rightleftarrows\C W\colon G$  is a weak symmetric monoidal Quillen adjunction  in the sense of \cite[Definition 3.6]{emmc}. Then the left adjoint $\mathcal LF$ in the derived adjunction \eqref{apderived}
is \emph{strong} symmetric monoidal, so it induces an adjoint pair
\begin{equation}\label{apderived2}
\xymatrix{\cat{\ho \C V}{}\ar@<.5ex>[r]^-{\mathcal{L}F}&\cat{\ho\C W}{}\ar@<.5ex>[l]^-{\mathcal{R}G}}
\end{equation}
given by the total derived functors in \eqref{apderived} on morphism objects.

\begin{prop}\label{conputa}
Suppose in addition that $\C V$ and $\C W$ satisfy the strong unit axiom 
and that $F$ satisfies the pseudo-cofibrant axiom and the $\unit$-cofibrant axiom. Consider
the following diagram of categories and functors,
$$\xymatrix{\cat{\C V}{}\ar[r]^-{F^{\catata}}\ar[d]_{p_{\C V}}&\cat{\C W}{}\ar[d]^{p_{\C W}}\\
\cat{\ho\C V}{}\ar[r]_{\mathcal LF}&\cat{\ho\C W}{}}$$
If we restrict to cofibrant objects in $\cat{\C V}{}$, then the square commutes up to natural isomorphism.
\end{prop}

Proposition \ref{putifar} follows directly from this one. Recall that the axioms in the statement of Proposition \ref{conputa} were introduced in \cite[Appendices A and B]{htnso2} in order to avoid cofibrancy hypotheses on tensor units.

\begin{proof}[Proof of Proposition \ref{conputa}]
Let $\CC H$ be a cofibrant $\C V$-category. The $\ho\C W$-category $\mathcal LFp_{\C V}\CC H$
has underlying graph $F(\CC H)$. Here we use that morphism objects in $\CC H$ are cofibrant
or $\unit$-cofibrant in the sense of \cite[Definition B.1]{htnso2}, see Remark \ref{comoson}, and that the left derived functor $\mathcal LF$ coincides with $F$ on ($\unit$-)cofibrant objects by the $\unit$-cofibrant axiom. The natural isomorphism $\mathcal LFp_{\C V}\CC H\cong p_{\C W} F^{\catata}$ is represented by $\chi_{\CC H}$ in \eqref{chi}, which is a weak equivalence in $\graph{\C W}{\ob\CC H}$ by Proposition \ref{SS}.
\end{proof}

There is a monoidal natural transformation
$$\xymatrix{\ho\C V\ar[rr]^{\mathcal LF}\ar[rd]_-{\pi_0}&&\ho\C W\ar[ld]^-{\pi_0}\\
&\operatorname{Set}&\ar@{=>}(15,-5);(20,-5)^-\mu}$$
defined by
$$\pi_0V=\ho\C V(\unit,V)\st{\mathcal LF}\To\ho\C W(\mathcal LF(\unit),\mathcal LF(V))\cong\ho\C W(\unit,\mathcal LF(V))=\pi_0\mathcal LF(V).$$
Obviously $\mu$ is a natural isomorphism if $\mathcal L F$ is fully faithful, e.g.~if $F\colon\C V\rightleftarrows\C W\colon G$ is a Quillen equivalence.

\begin{cor}\label{preserba}
Under the assumptions of Proposition \ref{conputa}, for any $\C V$-category $\CC H$ there is a natural functor $\pi_0\CC H\r\pi_0\mathcal L F^{\catata}(\CC H)$ which is the identity on objects. In particular $\mathcal L F^{\catata}$ sends $\C V$-intervals to $\C W$-intervals. The previous natural functor is an isomorphism if $F\colon\C V\rightleftarrows\C W\colon G$ is a Quillen equivalence. 
\end{cor}

We are now ready to prove the transfer theorem stated in the introduction.

\begin{proof}[Proof of Theorem \ref{transfer}]
Let us tackle the first part of the statement. It is enough to show that $F^{\catata}$ sends generating (trivial) cofibrations to (trivial) cofibrations. Since $F_S^{\catata}\dashv G_S$ is a Quillen
pair and $F^{\catata}$ is fiberwise defined by $F_S^{\catata}$, $F^{\catata}$ sends $T_{\{0,1\}}(I_{01})$ (resp.~$T_{\{0,1\}}(J_{01})$) to (trivial) cofibrations in $\cat{\C W}{}$, actually in $\cat{\C W}{\{0,1\}}$. The morphism $F^{\catata}(\varnothing\r \unit)$ coincides with the cofibration $\varnothing\r \unit$ in $\cat{\C W}{}$ since the left adjoints  $F^{\catata}$ and $F^{\catata}_{\{0\}}$ preserve initial objects.
Therefore, it is only left to check that $F^{\catata}(\theta_{\CC I})$ is a trivial cofibration for all generating $\C V$-intervals $\CC I$. 

What we have already showed, proves that $F^{\catata}$ preserves cofibrations, hence it is enough to check that $F^{\catata}(\theta_{\CC I})$  is a DK-equivalence. It is homotopically essentially surjective since $F^{\catata}(\tilde{\CC I})$ is a $\C W$-interval, see Proposition \ref{conputa} and Corollary \ref{preserba}. Let us see that it is also homotopically fully faithful. 
There is a commutative diagram in $\graph{\C W}{}$
$$\xymatrix{
F^{\catata}(\widetilde{i^*\CC I})\ar[r]^-{F^{\catata}(\theta_{\CC I})}&F^{\catata}(\tilde{\CC I})\\
F (\widetilde{i^*\CC I})\ar[u]^{\chi_{\widetilde{i^*\CC I}}}\ar[r]_-{F(\theta_{\CC I})}&F (\tilde{\CC I})\ar[u]_{\chi_{\tilde{\CC I}}}
}$$
The source and target of $\theta_{\CC I}$ are cofibrant in their respective fibers, hence the vertical arrows are weak equivalences by Proposition \ref{SS}. Since $\theta_{\CC I}$ is a cofibration with cofibrant source, it is a local cofibration by Corollary \ref{008} and Theorem  \ref{ict2}. Moreover, since $\theta_{\CC I}$ is a DK-equivalence it is actually a local trivial cofibration. Therefore $F(\theta_{\CC I})$ is also a local trivial cofibration, since $F$ is a left Quillen functor. Now  $F^{\catata}(\theta_{\CC I})$ is homotopically fully faithful by the commutativity of the previous diagram and the 2-out-of-3 property of weak equivalences in $\C V$.

Suppose now that $F\colon\C V\rightleftarrows\C W\colon G$ is a Quillen equivalence. Let $\varphi\colon \CC H\r G(\CC K)$ be a  $\C V$-functor with $\CC H$ cofibrant in $\cat{\C V}{}$ and $\CC K$ fibrant in $\cat{\C W}{}$. Denote by $\varphi'\colon F^{\catata}(\CC H)\r\CC K$ the adjoint $\C W$-functor of $\varphi$. The $\C V$-functors $\varphi$ and $\varphi'$ factor as
$$\CC H\st{\bar\varphi}\To\varphi^*G(\CC K)\st{c}\To G(\CC K),\qquad
F^{\catata}(\CC H) \st{\bar\varphi'}\To\varphi^*\CC K\st{c'}\To \CC K.$$
Here $c$ and $c'$ are cartesian, actually $G(c')=c$ since $G$ is a cartesian functor. Moreover, $\CC H$ (resp.~$\varphi^*\CC K$) is cofibrant (resp.~fibrant) in $\cat{\C V}{S}$ (resp.~$\cat{\C W}{S}$), $S=\ob\CC H$, see Remark \ref{comoson}. 
The $\C V$-functor $\varphi$ (resp.~$\varphi'$) is a DK-equivalence if and only if $\bar\varphi$ (resp.~$\bar\varphi'$) is a weak equivalence in $\cat{\C V}{S}$ (resp.~$\cat{\C W}{S}$) and $c$ (resp.~$c'$) is homotopically essentially surjective. Since $F_S^{\catata}\dashv G_S$ is a Quillen equivalence by Proposition \ref{qeq}, $\bar\varphi$ is a weak equivalence in $\cat{\C V}{S}$ if and only if $\bar\varphi'$ is a weak equivalence in $\cat{\C W}{S}$. Moreover, 
$\pi_0c$ is isomorphic to 
$\pi_0c'$ by Corollary \ref{preserba0}. This concludes the proof.
\end{proof}

\section{A kind of left properness result}

In this section we prove the following generalization of Theorem \ref{lp}, which shows that $\cat{\C V}{}$ shares some features with left proper model categories.

\begin{thm}\label{lp2}
If  $\C V$ satisfies the strong unit axiom \cite[Definition A.9]{htnso2},
$$\xymatrix{
\CC H\ar@{>->}[r]^\varphi\ar[d]^\sim_\phi\ar@{}[rd]|{\text{push}}&\CC K\ar[d]^{\phi'}\\
\CC H'\ar@{>->}[r]_{\varphi'}&\CC K'
}$$
is a push-out in  $\cat{\C V}{}$, 
$\CC H$ and $\CC H'$ are locally pseudo-cofibrant, $\varphi$ is a cofibration, and $\phi$ is a DK-equivalence, then $\phi'$ is also a DK-equivalence.
\end{thm}

\begin{proof}
The morphism $\phi'$ is homotopically essentially surjective by Proposition \ref{closure1}. We now  concentrate in showing that $\phi'$ is homotopically fully faithful. It is enough to assume that $\varphi$ is a push-out of a generating cofibration. This follows from \cite[Lemma A.17]{htnso2} and Proposition \ref{paralpc}.

The easy case is when $\varphi$ is a push-out of the generating cofibration is $\varnothing\r\unit$. In this case $\CC K$ is obtained from $\CC H$ by adding a disconnected object with trivial endomorphism monoid, and  $\CC K'$ is obtained from $\CC H'$  in the same way, so the result is obvious. The difficult case is when $\varphi$ is a push-out as in \eqref{pussy} with $f\in I$. In this other case, $\varphi$ and $\varphi'$ are the identity on objects and $\phi$ and $\phi'$ induce the same map $f\colon S\r S'$ on object sets. Moreover, the explicit description of such push-outs in Section \ref{pushpush} shows that the square in the statement decomposes as
$$\xymatrix{
\CC H\ar@{>->}[r]^\varphi\ar[d]^\sim_{\bar\phi}\ar@{}[rd]|{\text{push}}&\CC K\ar[d]^{\bar\phi'}\\
f^*\CC H'\ar@{>->}[r]_{f^*\varphi'}\ar[d]_{\text{cartesian}}&f^*\CC K'\ar[d]^{\text{cartesian}}\\
\CC H'\ar@{>->}[r]_{\varphi'}&\CC K'
}$$
Cartesian morphisms are fully faithful, hence it is enough to prove that $\bar\phi'$ is homotopically fully faithful. The upper square is in $\cat{\C V}{S}$, $\varphi$ is a cofibration, and $\bar\phi$ is a weak equivalence with locally pseudo-cofibrant source and target. 

Pseudo-cofibrant objecs are defined in arbitrary biclosed monoidal model categories such as $\graph{\C V}{S}$, see \cite[Definition A.1]{htnso2}. A $\C V$-graph is pseudo-cofibrant  if and only if it is locally pseudo-cofibrant. This can be easily checked by applying the criterion in \cite[Lemma A.3]{htnso2} to the sets of generating (trivial) cofibrations $I_S$ and $J_S$. 
The strong unit axiom for $\graph{\C V}{S}$ follows directly from the strong unit axiom for $\C V$. Then $\bar\phi'$ is a weak equivalence by \cite[Theorem D.13]{htnso2} for $\C C=\graph{\C V}{\ob\CC H}$ and $\mathcal O$ the associative operad, compare \cite[\S10]{htnso}.

\end{proof}


\begin{thebibliography}{MMSS01}

\bibitem[Amr11]{msctc}
I.~Amrani, \emph{{A model structure on the category of topological
  categories}}, \texttt{arXiv:1110.2695 [math.AT]}, October 2011.

\bibitem[Bat98]{hcctaimc}
M.~A. Batanin, \emph{{Homotopy coherent category theory and
  {$A_\infty$}-structures in monoidal categories}}, J. Pure Appl. Algebra
  \textbf{123} (1998), no.~1-3, 67--103. \MR{1492896 (99c:18007)}

\bibitem[Ber07]{mcssc}
J.~E. Bergner, \emph{{A model category structure on the category of simplicial
  categories}}, Trans. Amer. Math. Soc. \textbf{359} (2007), no.~5, 2043--2058.

\bibitem[BM13]{htec}
C.~Berger and I.~Moerdijk, \emph{{On the homotopy theory of enriched
  categories}}, Q. J. Math. \textbf{64} (2013), no.~3, 805--846.

\bibitem[Bor94]{borceux2}
F.~Borceux, \emph{{Handbook of categorical algebra 2}}, {Encyclopedia of Math.
  and its Applications}, no.~51, Cambridge University Press, 1994.

\bibitem[DK80]{slc}
W.~G. Dwyer and D.~M. Kan, \emph{{Simplicial localizations of categories}}, J.
  Pure Appl. Algebra \textbf{17} (1980), no.~3, 267--284.

\bibitem[FL81]{hin}
R.~Fritsch and D.~M. Latch, \emph{{Homotopy inverses for nerve}}, Math. Z.
  \textbf{177} (1981), no.~2, 147--179.

\bibitem[Gir71]{cna}
J.~Giraud, \emph{{Cohomologie non ab{\'e}lienne}}, {Die Grundlehren der
  mathematischen Wissenschaften}, vol. 179, Springer-Verlag, Berlin, 1971.

\bibitem[Hau13]{reic}
R.~Haugseng, \emph{{Rectification of enriched $\infty$-categories}}, \texttt{
  arXiv:1312.3881 [math.AT]}, December 2013.

\bibitem[Hir03]{hirschhorn}
P.~S. Hirschhorn, \emph{{Model categories and their localizations}},
  {Mathematical Surveys and Monographs}, vol.~99, American Mathematical
  Society, Providence, RI, 2003.

\bibitem[Hov99]{hmc}
M.~Hovey, \emph{{Model categories}}, {Mathematical Surveys and Monographs},
  vol.~63, American Mathematical Society, Providence, RI, 1999.

\bibitem[Hov13]{brew}
\bysame, \emph{{Brown representability and the Eilenberg-Watts theorem in
  homological algebra}}, 2013.

\bibitem[HSS00]{se}
M.~Hovey, B.~Shipley, and J.~Smith, \emph{{Symmetric spectra}}, J. Amer. Math.
  Soc. \textbf{13} (2000), no.~1, 149--208.

\bibitem[Joh02]{elephant1}
P.~T. Johnstone, \emph{{Sketches of an elephant: a topos theory compendium.
  {V}ol. 1}}, {Oxford Logic Guides}, vol.~43, The Clarendon Press Oxford
  University Press, New York, 2002.

\bibitem[Joy08]{tqca}
A.~Joyal, \emph{{The Theory of Quasi-Categories and its Applications}},
  \texttt{http://mat.uab.cat/~kock/crm/hocat/advanced-course/Quadern45-2.pdf},
  2008.

\bibitem[Kel05]{bcect}
G.~M. Kelly, \emph{{Basic concepts of enriched category theory}}, Repr. Theory
  Appl. Categ. (2005), no.~10, vi+137, Reprint of the 1982 original [Cambridge
  Univ. Press].

\bibitem[KL01]{vcatlplb}
G.~M. Kelly and S.~Lack, \emph{{{$\mathcal{V}$}-{C}at is locally presentable or
  locally bounded if {$\mathcal{V}$} is so}}, Theory Appl. Categ. \textbf{8}
  (2001), 555--575.

\bibitem[Lur09]{htt}
J.~Lurie, \emph{{Higher topos theory}}, {Annals of Mathematics Studies}, vol.
  170, Princeton University Press, Princeton, NJ, 2009.

\bibitem[MMSS01]{mcds}
M.~A. Mandell, J.~P. May, S.~Schwede, and B.~Shipley, \emph{{Model categories
  of diagram spectra}}, Proc. London Math. Soc. (3) \textbf{82} (2001), no.~2,
  441--512.

\bibitem[Mur11]{htnso}
F.~Muro, \emph{{Homotopy theory of nonsymmetric operads}}, Algebr. Geom. Topol.
  \textbf{11} (2011), 1541--1599.

\bibitem[Mur14]{htnso2}
\bysame, \emph{{Homotopy theory of non-symmetric operads, {II}: {C}hange of
  base category and left properness}}, Algebr. Geom. Topol. \textbf{14} (2014),
  229--281.

\bibitem[Ros09]{gbrhc}
J.~Rosick{\'y}, \emph{{Generalized {B}rown representability in homotopy
  categories}}, corrected version of Theory Appl. Categ. 14 (2005), 451--479,
  available at \texttt{arXiv:math/0506168v3 [math.CT] }, December 2009.

\bibitem[Shi07]{hzas}
B.~Shipley, \emph{{{$H\mathbb{Z}$}-algebra spectra are differential graded
  algebras}}, Amer. J. Math. \textbf{129} (2007), no.~2, 351--379.

\bibitem[SS00]{ammmc}
S.~Schwede and B.~Shipley, \emph{{Algebras and modules in monoidal model
  categories}}, Proc. London Math. Soc. (3) \textbf{80} (2000), no.~2,
  491--511.

\bibitem[SS03]{emmc}
\bysame, \emph{{Equivalences of monoidal model categories}}, Algebr. Geom.
  Topol. \textbf{3} (2003), 287--334 (electronic).

\bibitem[Sta13]{cmcpfo}
A.~E. Stanculescu, \emph{{Constructing model categories with prescribed fibrant
  objects}}, \texttt{arXiv:1208.6005 [math.CT]}, July 2013.

\bibitem[Tab05]{cmqdgcat}
G.~Tabuada, \emph{{Une structure de cat{\'e}gorie de mod{\`e}les de {Q}uillen
  sur la cat{\'e}gorie des dg-cat{\'e}gories}}, C. R. Math. Acad. Sci. Paris
  \textbf{340} (2005), no.~1, 15--19.

\bibitem[Tab09]{htsc}
\bysame, \emph{{Homotopy theory of spectral categories}}, Adv. Math.
  \textbf{221} (2009), no.~4, 1122--1143.

\bibitem[Tab10a]{dgvsc}
\bysame, \emph{{Differential graded versus simplicial categories}}, Topology
  Appl. \textbf{157} (2010), no.~3, 563--593.

\bibitem[Tab10b]{gscthhtm}
\bysame, \emph{{Generalized spectral categories, topological {H}ochschild
  homology and trace maps}}, Algebr. Geom. Topol. \textbf{10} (2010), no.~1,
  137--213.

\bibitem[TV08]{hagII}
B.~To{\"e}n and G.~Vezzosi, \emph{{Homotopical algebraic geometry. {II}.
  {G}eometric stacks and applications}}, Mem. Amer. Math. Soc. \textbf{193}
  (2008), no.~902, x+224.

\end{thebibliography}
 
 \providecommand{\bysame}{\leavevmode\hbox to3em{\hrulefill}\thinspace}
\providecommand{\MR}{\relax\ifhmode\unskip\space\fi MR }
\providecommand{\MRhref}[2]{%
  \href{http://www.ams.org/mathscinet-getitem?mr=#1}{#2}
}
\providecommand{\href}[2]{#2}

\end{document}